\documentclass[a4paper]{article}

\pdfoutput=1
\def\zibreport{0}

\usepackage{a4wide}

\usepackage[breaklinks,colorlinks=true,citecolor=black,linkcolor=black,urlcolor=black,hyperfootnotes=false]{hyperref}
\usepackage{url}

\usepackage{amsmath,amsfonts,amsthm,amssymb}
\usepackage[utf8]{inputenc}
\usepackage[numbers]{natbib}
\usepackage{doi}
\newcommand{\doilink}[1]{\href{http://dx.doi.org/#1}{\nolinkurl{doi:#1}}}

\usepackage{microtype}
\usepackage{booktabs}
\usepackage{longtable}
\usepackage{ifthen}

\usepackage{tikz}
\usepackage{pgfplots}
\pgfplotsset{compat=1.17}  %

\newcommand{\low}[1]{\underline{#1}}
\newcommand{\upp}[1]{\overline{#1}}
\newcommand{\Rbb}{\mathbb{R}}
\newcommand{\Rinf}{\overline{\Rbb}}
\newcommand{\Zbb}{\mathbb{Z}}
\newcommand{\abs}[1]{\lvert{#1}\rvert}

\title{Global Optimization of Mixed-Integer Nonlinear Programs with SCIP 8}

\author{
Ksenia Bestuzheva\footnote{
  Zuse Institute Berlin, Department AIS$^2$T, \texttt{bestuzheva@zib.de}, ORCID: 0000-0002-7018-7099},
Antonia Chmiela\footnote{
  Zuse Institute Berlin, Department AIS$^2$T, \texttt{chmiela@zib.de}, ORCID: 0000-0002-4809-2958},
Benjamin M\"uller\footnote{
  Zuse Institute Berlin, Department AIS$^2$T, \texttt{benjamin.mueller@zib.de}, ORCID: 0000-0002-4463-2873},\\
Felipe Serrano\footnote{
  Zuse Institute Berlin, Department AIS$^2$T, \texttt{serrano@zib.de}, ORCID: 0000-0002-7892-3951},
Stefan Vigerske\footnote{
  GAMS Software GmbH, c/o Zuse Institute Berlin, Department AIS$^2$T, \texttt{svigerske@gams.com}},
Fabian Wegscheider\footnote{
  Zuse Institute Berlin, Department AIS$^2$T}
}

\ifthenelse{\zibreport = 1}{
  \let\pdfoutorg\pdfoutput
  \let\pdfoutput\undefined
  \usepackage{zibtitlepage}
  \let\pdfoutput\pdfoutorg
  \ZTPTitle{Global Optimization of Mixed-Integer Nonlinear Programs with SCIP 8}
  \ZTPAuthor{Ksenia Bestuzheva, Antonia Chmiela, Benjamin M\"uller, Felipe Serrano, Stefan Vigerske, Fabian Wegscheider}
  \ZTPPreprint
  \ZTPNumber{22-XX}
  \ZTPMonth{12}
  \ZTPYear{2022}
  \ZTPInfo{}
}{}

\begin{document}

\ifthenelse{\zibreport = 1}{\zibtitlepage}{}

\maketitle
\abstract{
For over ten years, the constraint integer programming framework SCIP has been extended by capabilities for the solution of convex and nonconvex mixed-integer nonlinear programs (MINLPs).
With the recently published version~8.0, these capabilities have been largely reworked and extended.
This paper discusses the motivations for recent changes and provides an overview of features that are particular to MINLP solving in SCIP.
Further, difficulties in benchmarking global MINLP solvers are discussed and a comparison with several state-of-the-art global MINLP solvers is provided.}
\normalsize

\section{Introduction}

Mixed-integer nonlinear programming (MINLP) concerns with the optimization of an objective function such that a finite set of linear or nonlinear constraints and integrality conditions is satisfied.
The generality of this problem class means that many real-world applications can be modeled as MINLPs~\cite{Floudas1995,GrossmannKravanja1997,Pinter2006,TrespalaciosGrossmann2014}, but also that software that can handle this class efficiently becomes extremely complex.
Solvers for MINLP~\cite{BussieckVigerske2010} are often built on top of or by combining solvers for mixed-integer linear programming (MIP) and solvers that find locally optimal solutions for nonlinear programs (NLP).
In fact, one of the first commercial MINLP solvers, SCICONIC~\cite{Beale1980}, extends a MIP solver by piecewise linear approximations of low dimensional nonlinear terms.
The first general purpose solver was DICOPT~\cite{KocisGrossmann1989}, which decomposes the solution of an MINLP into a sequence of MIP and NLP solves~\cite{DuranGrossmann1986}, thereby building on established software for these two program classes.
DICOPT can solve MINLPs where nonlinear constraints are convex to optimality, but works only as a heuristic on nonconvex MINLPs.
The first general purpose solvers to solve also nonconvex MINLPs to optimality were $\alpha$BB, BARON, and GLOP~\cite{AdjimanFloudas1996,Sahinidis1996,SmithPantelides1999}, all based on convexification techniques for nonconvex constraints.
Also the solver SCIP (Solving Constraint Integer Programs), for which this paper provides an overview, belongs to the latter category.

In the following, MINLPs of the form
  \begin{equation}
    \begin{aligned}
      \min \quad& c^{\top}x, \\
      \mathrm{such\ that} \quad& \low{g} \leq g(x) \leq \upp{g}, \\
      & \low{b} \leq Ax \leq \upp{b}, \\
      & \low{x} \leq x \leq \upp{x}, \\
      & x_{\mathcal{I}} \in \Zbb^{\vert\mathcal{I}\vert},
    \end{aligned}
    \label{eq:minlp} \tag{MINLP}
  \end{equation}
are considered, where $\low{x}$, $\upp{x} \in \Rinf^{n}$, $\Rinf := \Rbb \cup \{\pm\infty\}$, $\low{x}\leq \upp{x}$, $\mathcal{I} \subseteq \{1, \ldots, n\}$, $c \in \Rbb^n$, $\low{g}$, $\upp{g}\in\Rinf^m$, $\low{g}\leq \upp{g}$,
$g : \Rbb^{n} \rightarrow \Rinf^m$ is specified explicitly in algebraic form,
$\low{b},\upp{b}\in\Rinf^{\tilde m}$, $\low{b}\leq\upp{b}$, and $A\in\mathbb{R}^{\tilde m\times n}$.
The restriction to a linear objective function is a technical detail of SCIP and without loss of generality.

The software SCIP has been designed as a branch-cut-and-price framework to solve different types of optimization problems, most generally \emph{constraint integer programs} (CIPs), and most importantly MIPs and MINLPs.
Roughly speaking, CIPs are finite-dimensional optimization problems with arbitrary constraints and a linear objective function that satisfy the following property: if all integer variables are fixed, the remaining subproblem must form a linear or nonlinear program.
The problem class of CIP was motivated by the modeling flexibility of constraint programming and the algorithmic requirements of integrating it with efficient solution techniques available for MIP~\cite{Achterberg2007}.

In order to solve CIPs, SCIP constructs relaxations -- typically linear programs (LPs).
If the relaxation solution is not feasible for the current subproblem, the plugins that handle the violated constraints need to take measures to eventually render the relaxation solution infeasible for the updated relaxation, for example by branching or separation~\cite{Achterberg2007}.
A plethora of additional plugin types, e.g., for presolving, finding feasible solutions, or tightening variable bounds, allow accelerating the solution process.
After 20 years of development of the framework itself and included plugins, SCIP includes mature solvers for MIP, MINLP, as well as several other problem classes.
Since November 2022, SCIP is freely available under an open-source license.

SCIP solves problems like \eqref{eq:minlp} to global optimality via a spatial branch-and-bound algorithm that mixes branch-and-infer and branch-and-cut~\cite{BeKiLeLiLuMa12}.
Important parts of the solution algorithm are presolving, domain propagation (that is, tightening of variable bounds), linear relaxation, and branching.
A distinguishing feature of SCIP is that its capabilities to handle nonlinear constraints are not limited to MINLPs, but can be used for any CIP.
For example, problems can be handled where linear and nonlinear constraints are mixed with typical constraints from constraint programming, as long as appropriate constraint handlers have been included in SCIP.
Since most constraint handlers in SCIP construct a linear relaxation of their constraints, also the handling of nonlinear constraints focuses on linear relaxations.
The emphasis on handling CIPs with nonlinear constraints rather than MINLP only is also a reason that the use of nonlinear relaxations or reformulations of complete MINLPs into other problem types, e.g., mixed-integer conic programs, has not been explored much so far.

The development of SCIP initially focused on solving CIPs where fixing all integer variables resulted in a linear program~\cite{Achterberg2007}.
However, it was soon realized that this requirement was not actually enforced by the implementation.
As long as constraint handlers were able to resolve infeasibilities by separation, branching, or other means, the problem could be handled by SCIP.
First experiments to handle nonlinear constraints in continuous variables were conducted for bilinear mixing constraints in mine production planning~\cite{BleyKochNiu2008}.
The positive results of these experiments motivated the decision to include support for more general nonlinear constraints.
With version 1.2~(2009), initial support for quadratic constraints (convex or nonconvex) and solving quadratically constrained programs (QCPs) to local optimality by Ipopt~\cite{WaechterBiegler2006} was added~\cite{BeHeVi09}.
For version 2.0~(2010), a primal heuristic that solves sub-MIPs was added~\cite{BertholdGleixner2014} and other large-neighborhood-search heuristics were extended to create sub-MINLPs~\cite{BeHePfVi11}.
Further, second-order cone constraints in three variables could be handled.
More general nonlinear constraints, specified in algebraic form, were first supported by SCIP~2.1~(2011)~\cite{Vigerske2013}.
Next to the specialized treatment for quadratic constraints, also handlers for signpower constraints ($x\abs{x}^p = z$ for some $p \geq 1$)~\cite{GlHeHuVi12} and 1-convex bivariate constraints ($f(x,y) = z$ for $f$ being convex or concave whenever $x$ or $y$ has been fixed)~\cite{BaMiVi12} were added.

With the basic handling of nonlinear constraints in place~\cite{VigerskeGleixner2016}, the next releases were dedicated to adding features that improved performance.
SCIP~3.0~(2012) brought optimization-based bound tightening (OBBT)~\cite{GleixnerBertholdMuellerWeltge2017} and an NLP diving heuristic.
SCIP~3.2~(2015) added a reformulation of general quadratic constraints into second-order cone constraints and separation for edge-concave decompositions of quadratic constraints~\cite{SCIPoptsuite32}.
With SCIP~4.0~(2017), higher-dimensional second-order cone constraints were disaggregated, KKT conditions for quadratic programs were utilized, multiple starting points were tried for NLP solves, solutions of the LP relaxation were projected onto a convex NLP relaxation, and also OBBT could be performed on the NLP instead of the LP relaxation~\cite{SCIPoptsuite40}.
Improved convexification of bilinear constraints by use of additional linear constraints~\cite{MuellerSerranoGleixner2020}, a new primal heuristic that solves a sequence of NLP reformulations, and interfaces to the NLP solvers filterSQP and Worhp~\cite{FletcherLeyffer1998,BueskensWassel2013,MuellerKuhlmannVigerske2017} were added for SCIP 5.0~(2017)~\cite{SCIPoptsuite50}.
The following two major releases brought a branch-and-price based solver for ring-packing~\cite{GleixnerMaherMuellerPedroso2020} (SCIP 6.0, 2018) and support for convex nonlinear subproblems in Benders Decomposition (SCIP 7.0~\cite{SCIPoptsuite70}, 2020).

That versions 6 and 7 added comparatively few features for MINLP was due to an ongoing complete overhaul on the way how nonlinear constraints were handled.
The primary motivation for this change, which was released with SCIP~8.0~(2022)~\cite{SCIPoptsuite80}, was to increase the reliability of the solver and to alleviate numerical issues that arose from problem reformulations and led to SCIP returning solutions that are feasible in the reformulated problem, but infeasible in the original problem.
More precisely, previous SCIP versions built an extended formulation of~\eqref{eq:minlp} explicitly, with the consequence that the original constraints were no longer included in the presolved problem.
Even though the formulations were theoretically equivalent, it was possible that $\varepsilon$-feasible solutions for the reformulated problem were not $\varepsilon$-feasible in the original problem.
SCIP~8 remedies this by building an implicit extended formulation as an annotation to the original problem.
A second motivation for the major changes in SCIP~8 was to reduce the ambiguity of expression and nonlinear structure types by implementing different plugin types for low-level structure types that define expressions, and high-level structure types that add functionality for particular, sometimes overlapping structures.
Finally, new features for improving the solver's performance on MINLPs were introduced with SCIP~8.
These include intersection, SDP (semi-definite programming), and RLT (reformulation linearization technique) cuts for quadratic expressions~\cite{ChmielaMunozSerrano2021,AchterbergBestuzhevaGleixner2022}, perspective strengthening~\cite{BestuzhevaGleixnerVigerske2021}, and symmetry detection~\cite{Wegscheider2019}.

SCIP can read MINLPs from files in the following formats: LP, MPS, NL (AMPL), OSiL, PIP, and ZIMPL.
In addition, problems can be passed to SCIP via interfaces to a variety of programming languages and modeling packages, including AMPL, C, GAMS, Java, Julia, Python, and MATLAB.

The following section provides an overview of the MINLP solving capabilities of SCIP.
Afterwards, the performance of SCIP is compared with that of other state-of-the-art global solvers for MINLP.

\section{MINLP capabilities of SCIP}

In the following, an overview of the facilities available in SCIP that are specific to the handling of MINLPs is provided.
First, available nonlinear functions are listed and the integration of nonlinear constraints into the branch-and-cut solver of SCIP is discussed.
Next, the concept of a \emph{nonlinear handler} is introduced, which is a new plug-in type that has been added with SCIP 8 and facilitates the integration of extensions that handle specific nonlinear structures.
The remainder of this section gives an overview of features available in SCIP that increase the efficiency of MINLP solving, e.g., cut generators to tighten the linear relaxation, presolve reductions to simplify the problem, and primal heuristics to find feasible solutions early.

To be concise, the presentation has been limited to high-level descriptions that spare technical details.
Unless specified otherwise, more details are often found in~\cite{SCIPoptsuite80}.

\subsection{Framework}

\subsubsection{Expressions}
\label{sec:expr}

\begin{figure}[ht]
  \centering
  \begin{tikzpicture}[scale=0.6, sibling distance=6em, level distance=6em]
    \tikzstyle{expr} += [shape=circle, draw, align=center, minimum size = 18pt] %

    \node[expr] (plus) {$\sum$}
      child[<-] {
        node[expr] (sqr1) {$\cdot^2$}
        child[<-] {
          node[expr] (log) {\scriptsize $\log$}
          child[<-] {
            node[expr] (x) {$x$}
          }
        }
      }
      child[<-] {
        node[expr] (prod) {$\prod$}
        edge from parent node[left,inner sep=.1em] {\scriptsize 2}
      }
      child[<-] {
        node[expr] (sqr2) {$\cdot^2$}
        child[<-] {
          node[expr] (y) {$y$}
        }
      };

    \draw[<-] (prod) to[out=-135,in=45] (log);
    \draw[<-] (prod) to[out=-45,in=135] (y);

    \node[yshift=+0.5cm] (t3) at (plus) {};

    \node[yshift=-0.5cm]               (t1) at (x) {};
    \node[yshift=-0.5cm]               (t2) at (y) {};
  \end{tikzpicture}
  \caption{Expression graph for algebraic expression $\log(x)^2 + 2\log(x)y+y^2$.}
  \label{fig:exprgraph}
\end{figure}
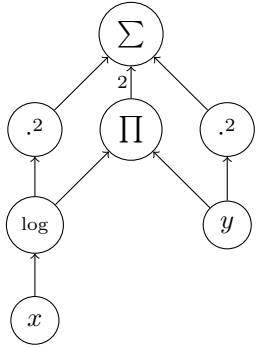

Algebraic expressions are well-formed combinations of constants, variables, and various algebraic operations such as addition, multiplication, and exponentiation, that are used to describe mathematical functions.
They are often represented by a directed acyclic graph with nodes representing variables, constants, and operations and arcs indicating the flow of computation, see Figure~\ref{fig:exprgraph} for an example.

Also in SCIP, expressions are stored as directed acyclic graphs, while all semantics of expression operands are defined by \textit{expression handler} plugins.
These handler implement callbacks that are used by methods in the SCIP core to manage expressions (create, modify, copy, free, parse, print), to evaluate and compute derivatives at a point, to evaluate over intervals, to simplify, to identify common subexpressions, to check curvature and integrality, and to iterate over it.
Some additional expression handler callbacks are used by the constraint handler for nonlinear constraints (Section~\ref{sec:consnl}) exclusively. %

Expression handlers for the following operators are included in SCIP~8.0:
\begin{itemize}
\item \texttt{val}: scalar constant;
\item \texttt{var}: a SCIP variable;
\item \texttt{sum}: an affine-linear function, $y\mapsto a_0 + \sum_{j=1}^k a_jy_j$ for $y\in\Rbb^k$ with constant coefficients $a\in\Rbb^{k+1}$;
\item \texttt{prod}: a product, $y\mapsto c\prod_{j=1}^ky_j$ for $y\in\Rbb^k$ with constant factor $c\in\Rbb$;
\item \texttt{pow}: a power with a constant exponent, $y\mapsto y^p$ for $y\in \Rbb$ and exponent $p\in\Rbb$ (if $p\not\in\Zbb$, then $y\geq 0$ is required);
\item \texttt{signpower}: a signed power, $y\mapsto \mathrm{sign}(y)\abs{y}^p$ for $y\in\Rbb$ and constant exponent $p\in\Rbb$, $p>1$;
\item \texttt{exp}: exponentiation, $y\mapsto \exp(y)$ for $y\in\Rbb$;
\item \texttt{log}: natural logarithm, $y\mapsto \log(y)$ for $y\in\Rbb_{>0}$;
\item \texttt{entropy}: entropy, $y\mapsto\begin{cases}-y\log(y), & \textrm{if }y > 0,\\0, & \textrm{if }y=0,\end{cases}$ for $y\in\Rbb_{\geq 0}$;
\item \texttt{sin}: sine, $y\mapsto\sin(y)$ for $y\in\Rbb$;
\item \texttt{cos}: cosine, $y\mapsto\cos(y)$ for $y\in\Rbb$;
\item \texttt{abs}: absolute value, $y\mapsto \abs{y}$ for $y\in\Rbb$.
\end{itemize}
In previous versions of SCIP, also high-level structures such as quadratic functions could be represented as expression types.
To avoid ambiguity and reduce complexity, this has been replaced by a recognition of quadratic expressions that is no longer made explicit by a change in the expression type.

\subsubsection{Constraint Handler for Nonlinear Constraints}
\label{sec:consnl}

All nonlinear constraints $\low{g}\leq g(x)\leq \upp{g}$ of \eqref{eq:minlp} are handled by the constraint handler for nonlinear constraints in SCIP, while the linear constraints $\low{b}\leq Ax\leq \upp{b}$ are handled by the constraint handlers for linear constraints and its specializations (e.g., knapsack, set-covering).
A constraint handler is responsible for checking whether solutions satisfy constraints and, if that is not the case, to resolve infeasibility by \emph{enforcing constraints}.
This applies in particular to solutions of the LP relaxation.
The nonlinear constraint handler currently enforces constraints by the following means:
\begin{description}
\item[DOMAINPROP] by analyzing the nonlinear constraints with respect to the variable bounds at the current node of the branch-and-bound tree, infeasibility or a bound tightening may be deduced, which allow pruning the node or cutting off the given solution, respectively; this is also known as \emph{domain propagation};
\item[SEPARATE] a cutting plane that is violated by the given solution may be computed;
\item[BRANCH] the current node of the branch-and-bound tree is subdivided, that is, a variable $x_i$ and a branching point $\tilde x_i\in[\low{x}_i,\upp{x}_i]$ are selected and two child nodes with $x_i$ restricted to $[\low{x}_i,\tilde{x}_i]$ and $[\tilde{x}_i,\upp{x}_i]$, respectively, are created.
\end{description}

To decide whether a node can be pruned (DOMAINPROP), an overestimate of the range of $g(x)$ with respect to current variable bounds is computed by means of interval arithmetics~\cite{Moore1966}.
If a constraint $k$ is found such that $g_k([\low{x},\upp{x}])\cap [\low{g}_k,\upp{g}_k]=\emptyset$, then there exists no point in $[\low{x},\upp{x}]$ for which this constraint is feasible.
A bound tightening may be computed by applying the same methods in reverse order.
That is, interval arithmetic is used to overestimate $g^{-1}([\low{g},\upp{g}])$, the preimage of $g(x)$ on $[\low{g},\upp{g}]$, and variable bounds are tightened to $[\low{x},\upp{x}]\cap g^{-1}([\low{g},\upp{g}])$.
This is also known as feasibility-based bound tightening (FBBT).
In the simplest case, callbacks of expression handlers are used to propagate intervals through expressions.
However, in some cases, other methods that take more structure into account or that use additional information to tighten variable bounds and constraint sides are used (see, e.g., Sections~\ref{sec:quadprop} and~\ref{sec:bilin}).

To construct a linear relaxation of the nonlinear constraints (SEPARATE option), an extended formulation is considered:
\begin{equation}
  \label{eq:minlp_ext}
  \tag{$\text{MINLP}_\text{ext}$}
  \begin{aligned}
    \min\; & c^\top x, \\
    \mathrm{such\ that}\; & h_i(x,w_{i+1},\ldots,w_{\hat m}) \lesseqgtr_i w_i, & i=1,\ldots,\hat m, \\
    & \low{b} \leq Ax \leq \upp{b}, \\
    & \low{x} \leq x \leq \upp{x}, \\
    & \low{w} \leq w \leq \upp{w}, \\
    & x_\mathcal{I} \in \Zbb^{\vert\mathcal{I}\vert}.
  \end{aligned}
\end{equation}
The functions $h_i$ are obtained from the expressions that define functions $g_i$ by recursively annotating subexpressions with auxiliary variables $w_{i+1},\ldots,w_{\hat m}$ for some $\hat m \geq m$.
Initially, slack variables $w_1,\ldots,w_m$ are introduced and assigned to the root of all expressions, i.e., $h_i:=g_i$, $\low{w}_i := \low{g}_i$, $\upp{w}_i:=\upp{g}_i$, for $i=1,\ldots,m$.
Next, for each function $h_i$, subexpressions $f$ may be assigned new auxiliary variables $w_{i'}$, $i'>m$, which results in extending~\eqref{eq:minlp_ext} by additional constraints $h_{i'}(x) = w_{i'}$ with $h_{i'} := f$.
Bounds $\low{w}_{i'}$ and $\upp{w}_{i'}$ are initialized to bounds on $h_{i'}$, if available.
Since auxiliary variables in a subexpression of $h_i$ always receive an index larger than $\max(m,i)$, the result is referred to by $h_i(x,w_{i+1},\ldots,w_{\hat m})$ for any $i=1,\ldots, \hat m$.
That is, to simplify notation, $w_{i+1}$ is used instead of $w_{\max(i,m)+1}$.
If a subexpression appears in several expressions, it is assigned at most one auxiliary variable and reindexing may be necessary to have $h_i$ depend on $x$ and $w_{i+1},\ldots, w_{\hat m}$ only.

For the (in)equality sense $\lesseqgtr_i$, a valid simplification would be to assume equality everywhere.
For performance reasons, though, it can be beneficial to relax certain equalities to inequalities if that does not change the feasible space of~\eqref{eq:minlp_ext} when projected onto $x$.
Therefore,
\[
  \lesseqgtr_i\; := \begin{cases}
  =,    & \text{if } \low{g}_i > -\infty,\; \upp{g}_i < \infty, \\
  \leq, & \text{if } \low{g}_i = -\infty,\; \upp{g}_i < \infty, \\
  \geq, & \text{if } \low{g}_i > -\infty,\; \upp{g}_i = \infty,
\end{cases}
\quad\text{ for }i=1,\ldots,m.
\]
For $i>m$, monotonicity of expressions is taken into account to derive $\lesseqgtr_i$. %

Whether to annotate a subexpression by an auxiliary variable depends on the structures that are recognized.
In the simplest case, every subexpression that is not already a variable is annotated with an auxiliary variable.
This essentially corresponds to the Smith Normal Form~\cite{SmithPantelides1999}.
For every function $h_i$ of~\eqref{eq:minlp_ext}, the callbacks of the corresponding expression handler can be used to compute linear under- and overestimators, such that a linear relaxation for~\eqref{eq:minlp_ext} is constructed.
It can, however, be beneficial to not add an auxiliary variable for every subexpression, thus allowing for more complex functions in~\eqref{eq:minlp_ext}.
This will be the discussed in Section~\ref{sec:nlhdlr} below.

\paragraph{Example}
Recall Figure~\ref{fig:exprgraph} and the constraint
\[
  \log(x)^2 + 2\log(x)\,y+y^2 \leq 4.
\]
By annotating the root of the expression graph with a slack variable $w_1$ and each other non-variable node with an auxiliary variable, the extended formulation
\begin{align*}
w_2 + 2w_3 + w_4 & \leq w_1, \\
w_5^2 & \leq w_2, \\
w_5\,y & \leq w_3, \\
y^2 & \leq w_4, \\
\log(x) & = w_5,  \\
w_1 & \leq 4.
\end{align*}
is obtained.
Bounds on auxiliary variables have been omitted here.
The constraints $w_5^2 = w_2$, $w_5y = w_3$, and $y^2 = w_4$ were relaxed to inequalities because $w_2+2w_3+w_4$ is monotonically increasing in each variable.
However, to relax $\log(x)=w_5$ to $\log(x)\leq w_5$, both $w_5^2$ and $w_5y$ would need to be monotonically increasing in $w_5$.
This would be the case if $\low{x}\geq 1$ and $\low{y}\geq 0$.
\medskip

If a constraint $h_i(x,w_{i+1},\ldots,w_{\hat m})\leq w_i$ (the $\geq$-case is analogous) of~\eqref{eq:minlp_ext} is violated and $h_i$ is nonconvex, then linear underestimators on $h_i$ can only be as tight as the convex envelope of $h_i$.
Therefore, it may not be possible find a hyperplane that is violated by the solution of the LP relaxation.
Since the convex envelope of $h_i$ depends on the bounds of variables appearing in $h_i$, these variables are candidates for branching (BRANCH).
More precisely, when an expression handler computes a linear under- or overestimator for $h_i(x,w_{i+1},\ldots,w_{\hat m})$, it also signals for which variables it used current variable bounds.
Marked original variables are then added to the list of branching candidates.
For an auxiliary variable $w_{i'}$, $i'>i$, the variables in the subexpression that $h_{i'}$ represents are considered for branching instead.

The decision on whether to add a cutting plane that separates the solution of the LP relaxation or to branch is rather complex, %
but the idea is to branch if either no cutting plane is found or if the violation of available cutting planes in the relaxation solution is rather small when compared to the convexification gap of the under/overestimators that define the cutting planes.
In the latter case, it may be beneficial to first reduce the convexification gap by branching.
To select one variable from the list of branching candidates, the violation of constraints in~\eqref{eq:minlp_ext} and historical information about the effect of branching on a given variable on the optimal value of the LP relaxation (``pseudo costs'') are taken into account. %
The branching point is a convex combination of the value of the variable in the LP relaxation and the mid-point of the variable's interval.

\subsubsection{Nonlinear Handlers}
\label{sec:nlhdlr}

In the previous example, four auxiliary variables were introduced to construct the extended formulation.
This is due to the expression handlers having a rather myopic view, basically, implementing techniques that can handle only their direct children.
It is clear that, for this example, an extended formulation that only replaces $\log(x)$ by an auxiliary variable~$w_2$ could be more efficient to solve.
However, this requires methods to detect the quadratic (or convex) structure and to either compute linear underestimators for the quadratic (convex) expression $w_2^2 + 2w_2y+y^2$ or to separate cutting planes for the set defined by $w_2^2 + 2w_2y+y^2\leq w_1$.

Such structure detection and handling methods are the task of the new \emph{nonlinear handler} plugins that were introduced with SCIP 8.
Nonlinear handlers determine the extended formulation~\eqref{eq:minlp_ext} by deciding when to annotate subexpressions with auxiliary variables.
That is, given a constraint $h_i(x) \lesseqgtr_i w_i$, a nonlinear handler analyses the expression that defines $h_i$ and attempts to detect specific structures.
At this point, it may also request to introduce additional auxiliary variables, thus changing $h_i(x)$ into $h_i(x,w_{i+1},\ldots,w_{\hat m})$.
In addition, it informs the constraint handler that it will now provide separation for $h_i(x,w_{i+1},\ldots,w_{\hat m}) \leq w_i$, or $\geq w_i$, or both.
If none of the nonlinear handlers declare that they will handle $h_i(x) \lesseqgtr_i w_i$, auxiliary variables are introduced for each argument of the root of the expression $h_i$ and expression handler callbacks are used to construct cutting planes from linear under-/overestimators.

In addition to separation, nonlinear handlers can also contribute to domain propagation.
This is implemented analogously to separation by setting up an additional extended formulation similarly to~\eqref{eq:minlp_ext}, with the main difference that slack and auxiliary variables are not actually created in SCIP and equalities are currently not relaxed to inequalities.

Note that the extended formulations are stored as \emph{annotation} on the original expressions.
Thus, for each task, the most suitable formulation can be used.
For example, feasibility is checked on the original constraints, domain propagation and separation use the corresponding extended formulations, but branching is performed, by default, with respect to original variables only.
With SCIP 7 and earlier, only one extended formulation was constructed explicitly and the connection to the original formulation was no longer available, leading to issues due to not ensuring that solutions are also ($\varepsilon$-)feasible for the original constraints.

In addition to the improved numeric reliability, the nonlinear handlers also allow for a higher flexibility when handling nonlinear structures.
For each node in an expression, more than one nonlinear handler can be attached, each one annotating possibly different subexpressions with auxiliary variables.
For example, for a nonconvex quadratic constraint $\sum_{i,j} a_{i,j} x_ix_j \leq w$, the nonlinear handler for quadratics can declare that it will provide separation (by intersection cuts, see Section~\ref{sec:intersection}), but that also other means of separation should be tried.
However, since no other nonlinear handler declares that it will provide separation, auxiliary variables are introduced for each argument of the sum, that is, an auxiliary variable $X_{ij}$ is assigned to each product $x_ix_j$.
For the corresponding constraints $x_ix_j\leq X_{ij}$ (if $a_{i,j}\geq 0$), the well-known McCormick underestimators~\cite{McCormick1976},
\begin{equation}
\label{eq:mccormick}
\begin{aligned}
X_{ij} & \geq \low{x}_ix_j + \low{x}_jx_i - \low{x}_i\low{x}_j, \\
X_{ij} & \geq \upp{x}_ix_j + \upp{x}_jx_i - \upp{x}_i\upp{x}_j,
\end{aligned}
\end{equation}
or other means (see Section~\ref{sec:bilin}) will be used to construct a linear relaxation.

\subsubsection{NLP Relaxation}
\label{sec:nlprelax}

Similar to the central LP relaxation of SCIP, an NLP relaxation is also available.
In contrast to constraint handlers, the NLP relaxation uses a common data structure to store its constraints.
At the moment, constraint handlers for linear constraints and the constraint handler for nonlinear constraints store a representation of their constraints in the NLP relaxation, so that in case of a MINLP, the NLP relaxation together with the integrality conditions on variables provides a unified view of the problem.
For nonlinear constraints, the original (non-extended) form $\low{g}\leq g(x)\leq \upp{g}$ is added to the NLP.
To find local optimal solutions for the NLP relaxation, interfaces to the NLP solvers filterSQP, Ipopt, and Worhp~\cite{FletcherLeyffer1998,WaechterBiegler2006,BueskensWassel2013} are available.
First- and second-order derivatives for these solvers are computed via CppAD~\cite{cppad}.

The NLP relaxation is mainly used by some primal heuristics~(Section~\ref{sec:heur}) and separators~(Section~\ref{sec:convextight}) at the moment.

\subsection{Presolving}
\label{sec:presolve}

When presolving nonlinear constraints, expressions are simplified and brought into a canonical form.
For example, recursive sums and products are flattened and fixed or aggregated variables are replaced by constants or sums of active variables.
In addition, it is ensured that if a subexpression appears several times (in the same or different constraints), always the same expression object is used.
This ensures that in the extended formulation \eqref{eq:minlp_ext} at most one auxiliary variable is attached to such common subexpressions.

\subsubsection{Variable Fixings}

Similar to what has been shown by Hansen et al.~\cite{HansenJaumardRuizXiong1993}, if a bounded variable $x_j$ does not appear in the objective ($c_j=0$), but in exactly one constraint $\low{g}_k \leq g_k(x) \leq \upp{g}_k$ where $g_k(x)$ is convex in $x_j$ for any fixing of other variables and $\upp{g}_k = +\infty$ (or concave in $x_j$ and $\low{g}_k=-\infty$), then there always exists an optimal solution where $x_j \in \{\low{x}_j,\upp{x}_j\}$.
For example, if $y\in[0,1]$ appears only in a constraint $xy+yz-y^2\leq 5$, then $y$ can be changed to a binary variable.

SCIP recognizes such variables for polynomial constraints (under additional assumptions~\cite{SCIPoptsuite80}) and changes the variable type to binary, if $\low{x}_j = 0$ and $\upp{x}_j=1$, or adds a bound disjunction constraint $x_j \leq \low{x}_j \vee x_j \geq \upp{x}_j$.
As a consequence, branching on $x_j$ leads to fixing the variable in both children.

\subsubsection{Linearization of Products}

The introduction emphasized that with SCIP 8, an explicit extended reformulation of nonlinear constraints is avoided.
An exception that proves this ``rule'' is the linearization of products of binary variables in presolving.
Doing so has the advantage that more of SCIP's techniques for MIP solving can be utilized.

In the simplest case, a product $\prod_i x_i$ is replaced by a new variable $z$ and a constraint of type ``and'' that models $z = \bigwedge_i x_i$ is added.
The ``and''-constraint handler will then separate a linearization of this product~\cite{BerHP09}.
For a product of only two binary variables, the linearization is added directly. %

For a quadratic function in binary variables with many terms, the number of variables introduced may be large.
Thus, in this case, a linearization that requires fewer additional variables is used, even though it may lead to a weaker relaxation. %

\subsubsection{KKT Strengthening for QPs}

A presolving method that aims to tighten the relaxation of a quadratic program~(QP) by adding redundant constraints derived from Karush-Kuhn-Tucker (KKT) conditions is available.
Consider a quadratic program of the form
\begin{align}
  \label{eq:QP}\tag{QP}
  \begin{aligned}
    \min\ \ &\tfrac{1}{2}\, x^\top Q x + c^\top x,\\
    \text{such that}\ \ & Ax \leq b,%
  \end{aligned}
\end{align}
where $Q \in \Rbb^{n \times n}$ is symmetric, $c \in \Rbb^n$, $A \in \Rbb^{m \times n}$, 
and $b \in \Rbb^m$. %
If~\eqref{eq:QP} is bounded, then all optima of~\eqref{eq:QP} satisfy the following KKT conditions:
\begin{align}
 \label{eq:KKT_QP}\tag{KKT}
 \begin{aligned}
    Q x + c + A^\top \mu & = 0,\\
    Ax & \leq b,\\
    \mu_i (Ax - b)_i & = 0, && \qquad i \in \{1, \dots, m\},\\
    \mu & \geq 0,
  \end{aligned}
\end{align}
where $\mu$ is the vector of Lagrangian multipliers of the constraints $Ax\leq b$.

In a specialized presolver, SCIP recognizes whether~\eqref{eq:minlp} is equivalent to~\eqref{eq:QP} by checking whether a quadratic objective function has been reformulated into a constraint.
If a~\eqref{eq:QP} has been found and all variables are bounded, then the equations~\eqref{eq:KKT_QP} are added as redundant constraints to the problem, whereby the complementarity constraints are formulated via special ordered sets of type 1.
The redundant constraints can help to strengthen the linear relaxation and prioritize branching decisions to satisfy the complementarity constraints, which focuses the search more on the local optima of~\eqref{eq:QP}.

In addition to~\eqref{eq:QP}, the implementation can also handle mixed-binary quadratic programs.
For all details, see~\cite{SCIPoptsuite40,Fischer2017}.
When this presolver was added to SCIP 4.0, it has shown to be very beneficial for box-constrained quadratic programs.
Due to the many changes and extensions in SCIP~8, in particular for the handling of quadratic constraints~(Section~\ref{sec:quad}), it needs to be reevaluated under which conditions this presolver should be enabled.
Currently, it is disabled by default.

\subsubsection{Symmetry Detection}
\label{sec:symmetry}

Symmetries in a MINLP are automorphisms on $\mathbb{R}^n$ that map optimal solutions to optimal solutions.
Such symmetries have an adverse effect on the performance of branch-and-bound solvers, because symmetric subproblems may be treated repeatedly.
Therefore, SCIP can enforce lexicographically maximal solutions from an orbit of symmetric solutions via bound tightening and separation of linear inequalities~\cite{HojnyPfetsch2019,SCIPoptsuite50,SCIPoptsuite70,SCIPoptsuite80}.

Since optimal solutions are naturally not known in advance, the symmetry detection resorts to find permutations of variables that map the feasible set onto itself and map each point to one with the same objective function value~\cite{Margot2010}.
These permutations are given by isomorphisms in an auxiliary symmetry detection graph, which is constructed from the problem data (e.g., $c$, $A$, $\mathcal{I}$, and the expressions that define $g(x)$)~\cite{Liberti2012a,Wegscheider2019}.

\subsection{Quadratics}
\label{sec:quad}

Since quadratic functions frequently appear in MINLPs (every second instance of MINLPLib~\cite{minlplib} has only linear and quadratic constraints), a number of techniques have been added to SCIP to handle this structure.
Next to the presolving methods that were discussed in the previous section, three nonlinear handlers and four separators deal with quadratic structures.
When none of the nonlinear handlers are active, then for each square and bilinear term in a quadratic function, an auxiliary variable is added in the extended formulation and gradient, secant, and McCormick under- and overestimators (see \eqref{eq:mccormick}) are generated.

\subsubsection{Domain Propagation}
\label{sec:quadprop}

If variables appear more than once in a quadratic function, then a term-wise domain propagation does not necessarily yield the best possible results, due to suffering from the so-called \textit{dependency problem} of interval arithmetics.
For example, it is easy to compute the range for $x^2+x$ for given bounds on $x$, or bounds on $x$ for a given interval on $x^2+x$, but standard interval arithmetics would treat the terms $x^2$ and $x$ separately, which can lead to overestimating the result.

Therefore, a specialized nonlinear handler in SCIP provides a domain propagation procedure for quadratics that aims to reduce overestimation. %
For this, the detection routine of the nonlinear handler writes a quadratic expression as
\begin{equation}
  \label{eq:quadexpr}
  q(y) = \sum_{i=1}^k q_i(y) \quad\text{with}\quad q_i(y) = a_iy_i^2 + c_iy_i + \sum_{j\in P_i} b_{i,j}y_iy_j,
\end{equation}
where $y_i$ is either an original variable ($x$) or another expression, $a_i,c_i\in\Rbb$, $b_{i,j}\in\Rbb\setminus\{0\}$, $j\in P_i \Rightarrow i\not\in P_j$ for all $j\in P_i$, $P_i\subset\{1,\ldots,k\}$, $i=1,\ldots,k$.
For functions $q_i$ with at least two terms (at least two of $a_i$, $b_{i,j}$, $j\in P_i$, and $c_i$ are nonzero), a relaxation is obtained by replacing each $y_j$ by $[\low{y}_j,\upp{y}_j]$, $j\in P_i$.
For this univariate quadratic interval-term in $y_i$, tight bounds can be computed~\cite{DomesNeumaier2010}.

In addition, bounds on variables $y_j$, $j\in P_i$, are computed by considering
\begin{equation}
 \label{eq:quadexpr2}
 \sum_{j\in P_i}b_{i,j}y_j \in ([\low{q},\upp{q}] - \sum_{i'\neq i} q_{i'}(y))/y_i - a_iy_i - c_i, \qquad y_i\in [\low{y}_i,\upp{y}_i],
\end{equation}
where $[\low{q},\upp{q}]$ are given bounds on $q(y)$.
After relaxing each $q_{i'}$ to an interval, bounds on the right-hand side of \eqref{eq:quadexpr2} are computed, which are then used to calculate bounds on each $y_j$, $j\in P_i$.

\subsubsection{Bilinear Terms}
\label{sec:bilin}

For a product $y_1y_2$, where $y_1$ and $y_2$ are either non-binary variables or other expressions, the expression handler for products already provides linear under- and overestimators and domain propagation that is best possible when considering the bounds $[\low{y}_1,\upp{y}_1] \times [\low{y}_2,\upp{y}_2]$ only.
However, if linear inequalities in $y_1$ and $y_2$ are available, then possibly tighter linear estimates and variable bounds can be computed.
In SCIP, this is done by a specialized nonlinear handler that implements the algorithm by Locatelli~\cite{Locatelli2018}.
The inequalities are found by projection of the LP relaxation onto variables $(y_1,y_2)$.
For more details, see~\cite{MuellerSerranoGleixner2020}.
An alternative method that uses linear constraints to tighten the relaxation of quadratic constraints are the RLT cuts described in the following.

\subsubsection{RLT Cuts}
\label{sec:rlt}

The Reformulation-Linearization Technique (RLT)~\cite{adams1986tight,adams1990linearization} has proven very useful to tighten relaxations of polynomial programming problems.
In SCIP, a separator of cuts that are computed via RLT for bilinear product relations in~\eqref{eq:minlp_ext} is available.

For simplicity, denote by $X_{ij}$ the auxiliary variable that is associated with a constraint $x_ix_j \lesseqgtr X_{ij}$ of \eqref{eq:minlp_ext} ($X_{ji}$ denotes the same variable as $X_{ij}$).
Recall that it is valid to replace $\lesseqgtr$ by $=$.
Given $X_{ij} = x_ix_j$, where $x_i \in [\low{x}_i,\upp{x}_i]$, $x_j \in [\low{x}_j,\upp{x}_j]$, and a linear constraint $a^\top x \leq b$, RLT cuts are derived by first multiplying the constraint by a nonnegative bound factors $(x_i - \low{x}_i)$,
$(\upp{x}_i - x_i)$, $(x_j - \low{x}_j)$, or $(\upp{x}_j - x_j)$.
For instance, consider multiplication by the factor $(x_i - \low{x}_i)$, which yields a valid nonlinear inequality:
\begin{equation}\label{eq:rlt_reformulated}
 a^\top x\, (x_i - \low{x}_i) \leq b\,(x_i - \low{x}_i).
\end{equation}
This is referred to as the reformulation step.

The linearization step is then performed for all terms $x_kx_i$ in~\eqref{eq:rlt_reformulated}.
If a product relation $X_{ki} = x_kx_i$ exists, then the product is replaced with $X_{ki}$.
If $x_k$ and $x_i$ are contained in the same clique, the product is replaced with an equivalent linear expression.
Otherwise, it is replaced by a linear under- or overestimator such as~\eqref{eq:mccormick}.

In addition, the RLT separator can reveal linearized products between binary and continuous variables.
To do so, it checks whether pairs of linear inequalities that are defined in the same triple of variables (one of them binary, the other two continuous) imply a product relation.
These implicit products can then be used in the linearization step of RLT cut generation~\cite{AchterbergBestuzhevaGleixner2022}.

\subsubsection{SDP Cuts}
\label{sec:sdp}

As in the previous section, denote by $X_{ij}$ the auxiliary variable that is associated with a constraint $x_ix_j \lesseqgtr X_{ij}$ of \eqref{eq:minlp_ext}.
A popular convex relaxation of the condition $X = xx^\top$ is given by requiring $X-xx^\top$ to be positive semidefinite.
Separation for the set $\{(x,X) : X-xx^\top \succeq 0\}$ itself is possible, but cuts are typically dense and may include variables $X_{ij}$ for products that do not exist in the problem. %
Therefore, only principal $2 \times 2$ minors of $X-xx^\top$, which also need to be positive semidefinite, are considered.
By Schur's complement, this means that the condition
\begin{equation}
   \label{eq:minorposdef}
	A_{ij}(x,X) := \begin{bmatrix} 1 & x_i & x_j \\ x_i & X_{ii} & X_{ij} \\ x_j & 
	X_{ij} & X_{jj} \end{bmatrix} \succeq 0
\end{equation}
needs to hold for any $i,j$, $i\neq j$.
A separator in SCIP detects minors for which $X_{ii}$, $X_{jj}$, $X_{ij}$ exist in \eqref{eq:minlp_ext} and enforces $A_{ij}(x,X)\succeq 0$.
To do so for a solution $(\hat{x},\hat{X})$ that violates~\eqref{eq:minorposdef}, an eigenvector $v\in\Rbb^3$ of $A_{ij}(\hat{x},\hat{X})$ with $v^\top A_{ij}(\hat{x},\hat{X})v<0$ is computed and the globally valid linear inequality $v^\top A_{ij}(x,X)v \geq 0$ is added.

\subsubsection{Intersection Cuts}
\label{sec:intersection}

Intersection cuts~\cite{Tuy1964,Balas1971} have shown to be an efficient tool to strengthen relaxations of MIPs.
Recently, Mu{\~{n}}oz and Serrano showed how to compute the tightest possible intersection cuts for quadratic programs~\cite{MunozSerrano2020}.
This method has been implemented in SCIP~\cite{ChmielaMunozSerrano2021}.

Assume a nonconvex quadratic constraint of~\eqref{eq:minlp_ext} is $q(y)\leq w$ with $q(y)$ as in \eqref{eq:quadexpr} and $w$ an auxiliary variable.
The separation of intersection cuts is implemented for the set $S := \{ (y,w) \in \Rbb^k : q(y) \leq w \}$ that is defined by this constraint.

Let $(\hat{y},\hat w)$ be a basic feasible LP solution violating $q(y) \leq w$.
First, a convex inequality $g(y,w) < 0$ is build that is satisfied by $(\hat{y},\hat w)$, but by no point of $S$.
This defines a so-called \emph{$S$-free set} $C = \{ (y,w) \in \Rbb^{k+1} : g(y,w) \leq 0 \}$, that is, a convex set with $(\hat{y},\hat w) \in \text{int}(C)$ containing no point of $S$ in its interior.
The quality of the resulting cut highly depends on which $S$-free set is used, but using \textit{maximal} $S$-free sets yield the tightest possible intersection cuts~\cite{MunozSerrano2020}.

\begin{figure}[ht]
\centering
\includegraphics[width=.3\linewidth]{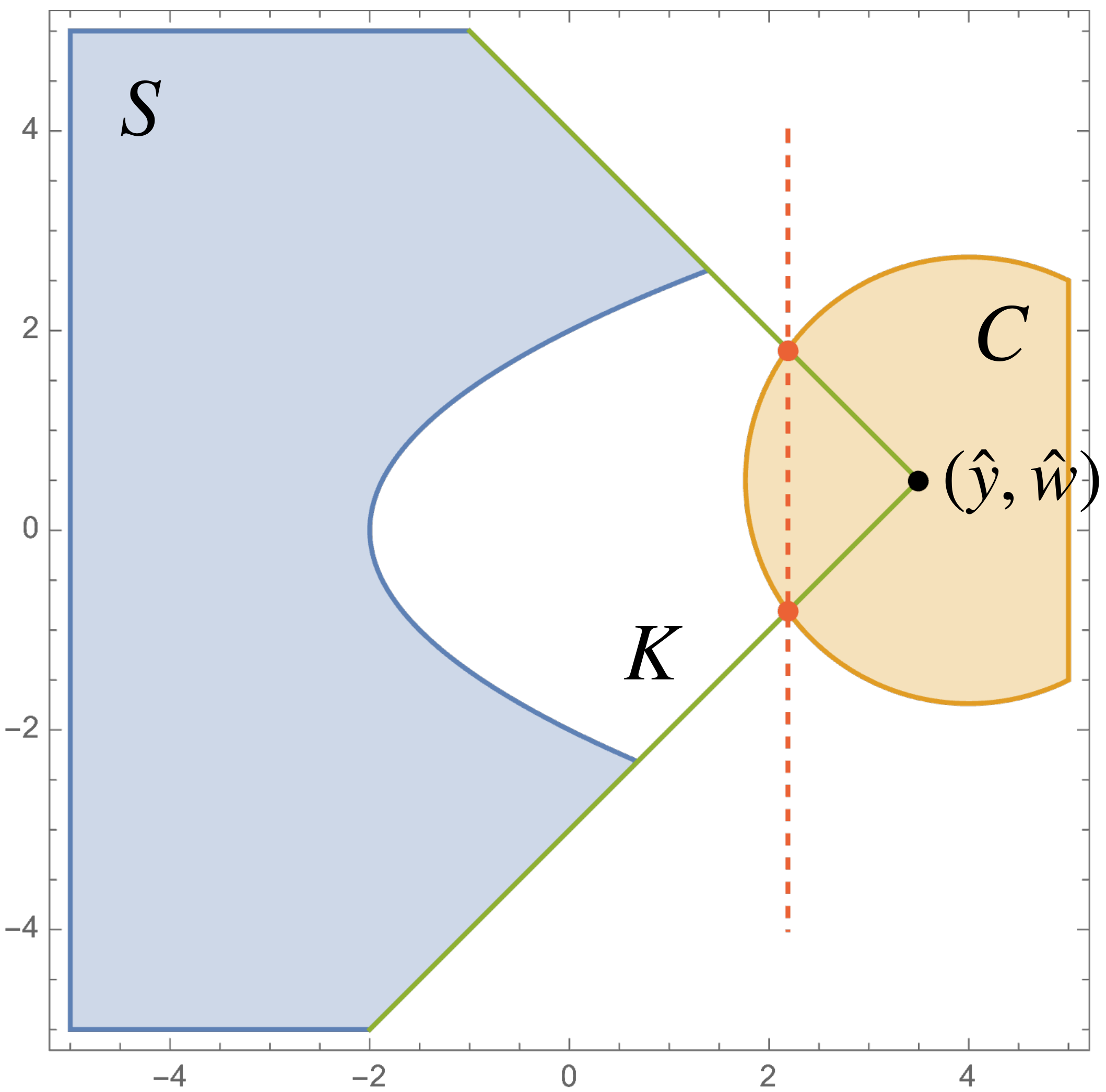}
\caption{An intersection cut (red) separating the basic feasible LP solution $(\hat y,\hat w)$ from $S$ (blue). The cut is computed using the intersection points of an $S$-free set $C$ (orange) and the rays of a simplicial cone $K \supseteq S$ (boundary in green) with apex $(\hat y,\hat w)\not\in S$.}
\label{fig:intercut}
\end{figure}

By using the conic relaxation $K$ of the LP-feasible region defined by the nonbasic variables at $(\hat{y},\hat w)$, the intersection points between the extreme rays of $K$ and the boundary of $C$ are computed.
The intersection cut is then defined by the hyperplane going through these points and successfully separates  $(\hat{x},\hat w)$ and $S$.
See Figure~\ref{fig:intercut} for an illustration.
To obtain even better cuts, there is also a strengthening procedure implemented that uses the idea of negative edge extension of the cone $K$~\cite{Glover1974}.

In addition to the separation of intersection cuts for a set $S$ given by a constraint $q(y)\leq w$, SCIP can also generate intersection cuts for implied quadratic equations.
Recall the matrix of auxiliary variables $X$ as introduced in Section~\ref{sec:rlt}.
The condition $X=xx^\top$ implies that $X$ needs to have rank 1.
Therefore, any $2\times 2$ minor $\begin{pmatrix}
X_{i_1j_1} & X_{i_1j_2} \\
X_{i_2j_1} & X_{i_2j_2}
\end{pmatrix}$
of $X$ needs to have determinant zero.
That is, for any set of variable indices $i_1$, $i_2$, $j_1$, $j_2$ with $i_1\neq i_2$ and $j_1\neq j_2$, the condition
$
  X_{i_1j_1}X_{i_2j_2} = X_{i_1j_2}X_{i_2j_1}
$
needs to hold.
If all variables in this condition exist in~\eqref{eq:minlp_ext} and a solution violates this condition, then the previously described procedure to generate intersection cuts is applied to the set defined by this condition.

Since intersection cuts can be rather dense, it is not clear yet how to decide when it will be beneficial to generate such cuts.
Their separation is therefore currently disabled by default.
For more details, see~\cite{ChmielaMunozSerrano2021}.

\subsubsection{Edge-Concave Cuts}

Another method to obtain a linear outer-approximation for a quadratic constraint is by utilizing an edge-concave decomposition of the quadratic function.
This has shown to be particularly useful for randomly generated quadratic instances~\cite{MisenerFloudas2012,MisenerSmadbeckFloudas2015}.
A function is edge-concave over the variables' domain (e.g., $[\low{x},\upp{x}]$) if it is componentwise concave.

Given a quadratic function, the separator for edge-concave cuts solves an auxiliary MIP to partition the square and bilinear terms into a sum of edge-concave functions and a remaining function.
Since the convex envelope of edge-concave functions is \emph{vertex-polyhedral}~\cite{Tardella2004}, that is, it is a polyhedral function with vertices corresponding to the vertices of the box of variable bounds, facets on the convex envelope of each edge-concave function can be computed by solving an auxiliary linear program (see also Section~\ref{sec:convexnlhdlr}).
For the function of remaining terms, term-wise linear underestimators such as~\eqref{eq:mccormick} are summed up.

Since the current implementation of edge-concave cuts in SCIP has not shown to be particularly useful for general MINLP, this separator is disabled for now.

\subsubsection{Second-Order Cones}

An important connection between MINLP and conic programming is the detection of constraints that can be represented as a second-order cone (SOC) constraint, since the latter defines a convex set, while the original constraint may use a nonconvex constraint function.

A specialized nonlinear handler aims to detect SOC representable structures.
In the detection phase, a constraint $h_i(x) \leq w_i$ (the case $\geq$ is handled similarly) of the extended formulation~\eqref{eq:minlp_ext} is passed to the nonlinear handler.
For this constraint, it is checked whether it defines a bound on an Euclidian norm ($\sqrt{\sum_{j=1}^k (a_j y_j^2 + b_j y_j) + c}\leq w_i$ for some coefficients $a_j,b_j,c\in\Rbb$, $a_j>0$, where $y_j$ is either an original variable or some subexpression of $h_i(\cdot)$), %
or is a quadratic constraint that is SOC-representable~\cite{MahajanMunson2010}.
Since the introduction of slack variables $w_i$, $i\leq m$, may prevent such a detection, the equivalent constraint $h_i(x) \leq \bar{w}_i$ is considered instead. %

A detected SOC constraint is stored in the form
\begin{equation}
 \label{eq:soc}
 \sqrt{\sum_{j=1}^k (v_j^\top y + \beta_j)^2} \leq v_{k+1}^\top y + \beta_{k+1}
\end{equation}
with $v_j\in\Rbb^\ell$, $j=1,\ldots,k+1$, where $y_1,\ldots,y_\ell$ are variables of~\eqref{eq:minlp_ext}.
Since the left-hand side of~\eqref{eq:soc} is convex, a solution $\hat y$ that violates~\eqref{eq:soc} can be separated by linearization of the left-hand side of~\eqref{eq:soc}.

However, if there are many terms on the left-hand side of~\eqref{eq:soc} ($k$ being large), then it can require many cuts to provide a tight linear relaxation of~\eqref{eq:soc}.
Thus, a disaggregation of the cone~\cite{Vielma2016} is used if $k\geq 3$:
\begin{align}
 (v_j^\top y + \beta_j)^2 & \leq z_j (v_{k+1}^\top y + \beta_{k+1}), \quad j=1,\ldots,k, \label{eq:socext1} \\
 \sum_{j=1}^k z_j & \leq v_{k+1}^\top y + \beta_{k+1}, \label{eq:socext2}
\end{align}
where variables $z_1,\ldots,z_k$ are new variables. %
A solution $(\hat y,\hat z)$ that violates~\eqref{eq:soc} needs to violate also~\eqref{eq:socext1} for some $j\in\{1,\ldots,k\}$ or \eqref{eq:socext2}.
The latter is already linear and can be added as a cut.
If a rotated second-order cone constraint~\eqref{eq:socext1} is violated for some $j$, then it is transformed into the standard form
\[
 \sqrt{4(v_j^\top y + \beta_j)^2 + (v_{k+1}^\top y + \beta_{k+1} - z_j)^2} \leq v_{k+1}^\top y + \beta_{k+1} + z_j
\]
and a gradient cut is constructed by linearization of the left-hand side.

\subsection{Convexity}
\label{sec:convex}

\subsubsection{Convex and Concave Constraints}
\label{sec:convexnlhdlr}

For the linear underestimation of functions like $x\exp(x)$ or $x^2 + 2xy + y^2$, the construction of an extended formulation ($xw$, $\exp(x)=w$; $w_1+2w_2+w_3$, $w_1=x^2$, $w_2=xy$, $w_3=y^2$) is not advisable.
Instead, hyperplanes that support the epigraph of a convex function can be used if convexity is recognized.
In SCIP, specialized nonlinear handlers are available to detect for a function $h_i(x)$ of~\eqref{eq:minlp_ext} the subexpressions that need to be replaced by auxiliary variables $w_{i+1},\ldots, w_{\hat{m}}$ such that the remaining expression $h_i(x,w_{i+1},\ldots, w_{\hat{m}})$ is convex or concave.
The detection utilizes the often applied rules for convexity/concavity of function compositions (e.g., $f$ convex and monotone decreasing, $g$ concave $\Rightarrow$ $f \circ g$ convex), but applies them in reverse order.
That is, instead of deciding whether a function is convex/concave based on information on the convexity/concavity and monotonicity of its arguments, the algorithm formulates conditions on the convexity/concavity of the function arguments given a convexity/concavity requirement on the function itself.
When a condition on an argument cannot be fulfilled, it is replaced by an auxiliary variable.

Next to ``myopic'' rules for convexity/concavity that are implemented by the expression handlers, also rules for product compositions ($af(b g(x)+c) g(x)$ with constants $a,b,c$ and repeating subexpression $g(x)$), signomials ($c\prod_{j=1}^k f_j^{p_j}(x)$ with $c,p_j\in\Rbb$ and subexpressions $f_j(x)$, $j=1,\ldots,k$), and quadratic forms are available.
The latter may check for definiteness of its Hessian by calculating its eigenvalues.
Further, it has been shown that for a composition of convex functions $f \circ g$, it can be beneficial for the linear relaxation to consider the extended formulation $f(w)$, $w\geq g(x)$, instead of the composition $f(g(x))$~\cite{TawarmalaniSahinidis2005}.
This is enforced by a small variation of the detection algorithm. %

When a convex constraint $h_i(x,w_{i+1},\ldots,w_{\hat m}) \leq w_i$ of~\eqref{eq:minlp_ext} is violated at a point $(\hat x,\hat w)$, a tangent on the graph of $h_i$ at $(\hat x,\hat w)$ is used to compute a separating hyperplane.
The slope of the tangent is given by the gradient of $h_i$ at $(\hat x,\hat w)$, which is calculated via automatic differentiation on the expression graph.
If, however, $h_i$ is univariate, that is, $h_i(x,w_{i+1},\ldots,w_{\hat m})=f(y)$ for some variable $y$, and $y$ is integral, then taking the hyperplane through the points $(\lfloor \hat y\rfloor, f(\lfloor \hat y\rfloor))$ and $(\lfloor \hat y\rfloor+1, f(\lfloor \hat y\rfloor+1))$ can give a tighter underestimator.

For a concave function $h_i(x,w_{i+1},\ldots,w_{\hat m})$, any hyperplane $\alpha x+\beta w+\gamma$ that underestimates $h_i(x,w_{i+1},\ldots,w_{\hat m})$ in all vertices of the box $[\low{x},\upp{x}]\times[\low{w}_{i+1},\upp{w}_{i+1}]\times\cdots\times[\low{w}_{\hat m},\upp{w}_{\hat m}]$ is a valid linear underestimator, since $h_i$ is vertex-polyhedral with respect to the box.
Maximizing $\alpha \hat x + \beta \hat w + \gamma$ such that $\alpha x + \beta w + \gamma$ does not exceed $h_i(x,w_{i+1},\ldots,w_{\hat m})$ for all vertices gives an underestimator that is as tight as possible at a given reference point $(\hat x, \hat w)$.
For the frequent cases $k=1$ and $k=2$, routines that directly compute such an underestimator are available.
For $k>2$, a linear program is solved.
Since the size of this LP is exponential in $k$, underestimators for concave functions in more than 14 variables are currently not computed.

\subsubsection{Tighter Gradient Cuts}
\label{sec:convextight}

The separating hyperplanes generated for convex functions of~\eqref{eq:minlp_ext} as discussed in the previous section are, in general, not supporting for the feasible region of~\eqref{eq:minlp_ext}, because the point where the functions are linearized is not at the boundary of the feasible region (which is the reason why it needs to be separated).
Therefore, often several rounds of cut generation and LP solving are required until the relaxation solution satisfies the convex constraints.
Solvers for convex MINLP have handled this problem in various ways~\cite{DuranGrossmann1986,KronqvistLundellWesterlund2016}, but the basic idea is to build gradient cuts at a suitable boundary point of the feasible region.

In SCIP, three procedures for building tighter and/or deeper gradient cuts for convex relaxations are included.
The first two methods compute a point on the boundary of the set defined by all convex constraints of~\eqref{eq:minlp} that is close to the point to be separated.
The first method solves an additional nonlinear program to project the point to be separated onto the convex set.
Since solving an NLP for every point to be separated can be quite expensive, the second method, going back to an idea by Veinott~\cite{Veinott1967}, does a binary search between an interior point of the convex set and the point to be separated.
The interior point is computed once in the beginning of the search by solving an auxiliary NLP.
For more details, see~\cite{SCIPoptsuite40}.

The third method does not aim to separate a given point, but utilizes the feasible points that are found by primal heuristics of SCIP.
When a new solution is found, gradient cuts are generated at this solution for convex constraints of~\eqref{eq:minlp_ext} and added to the cutpool.
If such a cut is later found to separate the relaxation solution, it is added to the LP.

All methods are currently disabled as they require more tuning to be efficient in general.

\subsection{Quotients}
\label{sec:quotient}

Note that the available expression handlers (see Section~\ref{sec:expr}) do not include a handler for quotients, since they can equivalently be written using a product and a power expression.
Therefore, the default extended formulation for an expression $y_1y_2^{-1}$ is given by replacing $y_2^{-1}$ by a new auxiliary variable $w$.
The linear outer-approximation is then obtained by estimating $y_1w$ and $y_2^{-1}$ separately.
However, tighter linear estimates are often possible.
Therefore, a specialized nonlinear handler checks whether a given function $h_i(x)$ can be cast as
\begin{equation}\label{eq:quotient_constraint}
 f(y) = \frac{ay_1 + b}{cy_2 + d} + e
\end{equation}
with $a,b,c,d,e\in\Rbb$, $a,c\neq 0$, and $y_1$ and $y_2$ being either original variables or subexpressions of $h_i(x)$.

Tight linear estimators for~\eqref{eq:quotient_constraint} are computed by distinguishing a number of cases.
For example, for $a\low{y}_1 +b\geq 0$ and $c\low{y}_2+d > 0$ (if $c>0$), a linear underestimator is obtained by computing a tangent on the graph of the convex underestimator of $f$ that is given by~\cite{ZamoraGrossmann1998}.
A linear overestimator is obtained by computing a facet on the concave envelope of $f$, which is easy since $-f$ is vertex-polyhedral. %
Furthermore, in the univariate case ($y_1=y_2$), $f$ is either convex or concave on $[\low{y}_1,\upp{y}_1]$ if $-d/c\not\in[\low{y}_2,\upp{y}_2]$.

Since in the univariate case the same variable appears twice, also a specialized domain propagation method that avoids the dependency problem of interval arithmetic is available.

\subsection{Perspective Strengthening}

Perspective reformulations have shown to efficiently tighten relaxations of convex mixed-integer nonlinear programs with on/off-structures, which are often modeled via big-M constraints or semi-continuous variables~\cite{frangioni2006perspective}.
A variable $x_j$ is semi-continuous with respect to the binary indicator variable $x_{j'}$, $j'\in \mathcal{I}$, if it is restricted to the domain $[\low{x}^1_j, \upp{x}^1_j]$ when $x_{j'}=1$ and has a fixed value $x^0_j$ when $x_{j'}=0$.

In SCIP, a strengthening of under- and overestimators for functions that depend on semi-continuous variables is available.
Consider a constraint $h_i(x,w_{i+1},\ldots,w_{\hat m}) \lesseqgtr w_i$ of~\eqref{eq:minlp_ext} and write $h_i$ as a sum of its nonlinear and linear parts:
\[
  h_i(x,w_{i+1},\ldots,w_{\hat m}) = h_i^\text{nl}(x_\text{nl},w_\text{nl}) + h_i^\text{l}(x_\text{l},w_\text{l}),
\]
where $h_i^\text{nl}$ is a nonlinear function, $h_i^\text{l}$ is a linear function, $x_\text{nl}$ and $w_\text{nl}$ are the vectors of variables $x$ and $w$, respectively, that appear only in the nonlinear part of $h_i$, and $x_\text{l}$ and $w_\text{l}$ are the vectors of variables $x$ and $w$, respectively, that appear only in the linear part of $h_i$.
A strengthening of under- or overestimators for $h_i(x,w_{i+1},\ldots,w_{\hat m})$ is attempted if $x_\text{nl}$ and $w_\text{nl}$ are semi-continuous with respect to the same indicator variable $x_{j'}$.

To determine whether a variable $x_j$ is semi-continuous, bounds on $x_j$ that are implied by fixing a binary variable are analyzed.
The implied bounds can be obtained either from linear constraints directly or by probing, and are stored by SCIP in a globally available data structure.
If a pair of implied bounds on $x_j$ with the same binary variable $x_{j'}$ is found, i.e.,
\begin{align*}
x_j &\leq \alpha^{(u)} x_{j'} + \beta^{(u)},\\
x_j &\geq \alpha^{(\ell)} x_{j'} + \beta^{(\ell)},
\end{align*}
and $\beta^{(u)} = \beta^{(\ell)}$, then $x_j$ is a semi-continuous variable with $x_j^0 = \beta^{(u)}$, $\low{x}^1_j = \alpha^{(\ell)} + \beta^{(\ell)}$, and $\upp{x}^1_j = \alpha^{(u)} + \beta^{(u)}$.
In addition, an auxiliary variable $w_i$ is found to be semi-continuous if function $h_i(x,w_{i+1},\ldots,w_{\hat m})$ depends only on semi-continuous variables with the same indicator variable.

\begin{figure}[ht]
\includegraphics[width=0.49\linewidth]{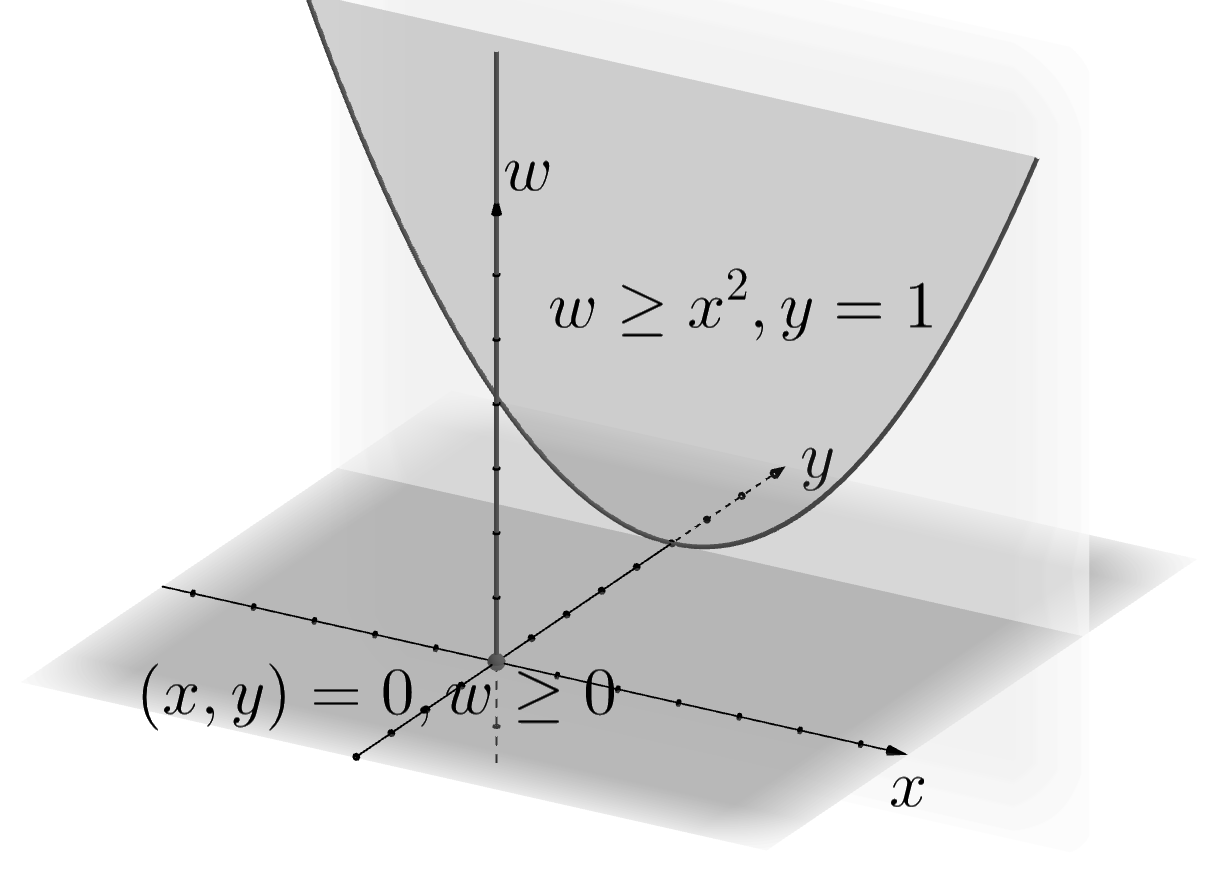}\hfill
\includegraphics[width=0.49\linewidth]{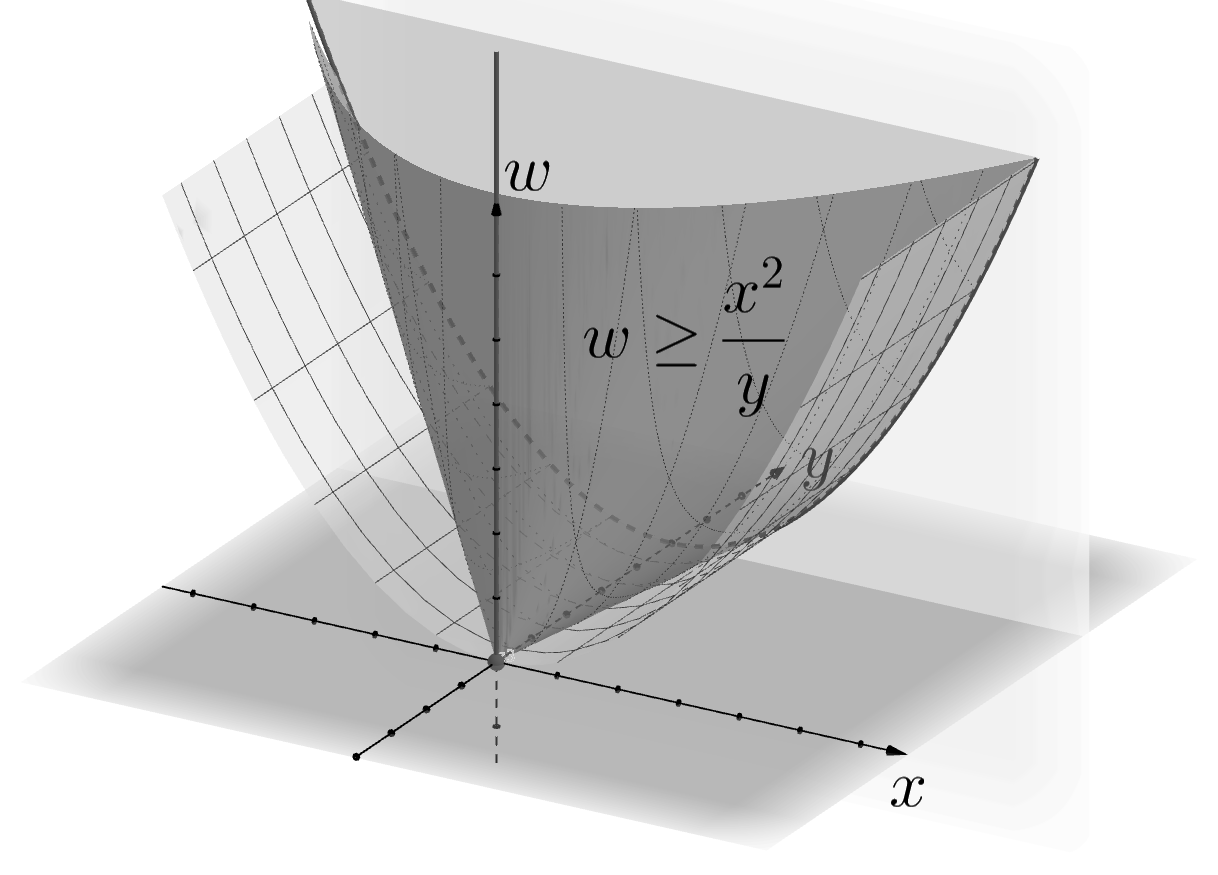}
\caption{The original set $\{(x,y,w): x^2 \leq w, y\in\{0,1\}, y = 0 \rightarrow x = 0\}$ (left) and a continuous relaxation given by $\{(x,y,w): x^2 \leq wy, y\in[0,1], w\geq 0\}$ (right).
From the original set, cuts of the form $\hat x^2 + 2\hat x (x-\hat x)\leq w$ for some reference point $\hat x$ would be generated.
With perspective strengthening, a linearization on the right set is obtained instead, i.e., $\hat x^2 + 2\hat x (x-\hat x) + \hat x^2 (1-y) \leq w$.
The latter is typically better as it is tight for $y=0$ as well.}
\label{fig:persp}
\end{figure}

Assume that a linear underestimator $\ell(x,w_{i+1},\ldots,w_{\hat m}$) of $h_i(x,w_{i+1},\ldots,w_{\hat m})$ has been computed and split it into parts corresponding to the nonlinear and linear variables of $h_i$, respectively:
\[
  \ell(x,w_{i+1},\ldots,w_{\hat m}) =
  \ell^\text{nl}(x_\text{nl},w_\text{nl}) + \ell^\text{l}(x_\text{l},w_\text{l}).
\]
The perspective strengthening consists of extending the part of the underestimator that corresponds to the nonlinear part such that it is tight for $x_{j'}=0$:%
\[
  \ell^\text{nl}(x_\text{nl},w_\text{nl}) +
  \left(h^\text{nl}_i(x^0_\text{nl},w^0_\text{nl}) - \ell^\text{nl}(x^0_\text{nl},w^0_\text{nl})\right)(1-x_{j'}) +
  \ell^\text{l}(x_\text{l},w_\text{l}).
\]
The linear part remains unchanged, since it shares none of the variables with the nonlinear part.
This extension ensures that the estimator is equal to $h_i(x,w_{i+1},\ldots,w_{\hat m})$ for $x_{j'}=0$, $x_\text{nl} = x^0_\text{nl}$, and $w_\text{nl} = w^0_\text{nl}$, and equal to $\ell(x,w_{i+1},\ldots,w_{\hat m})$ for $x_{j'}=1$.
If $h_i$ is convex, cuts obtained this way are equivalent to the classic perspective cuts~\cite{frangioni2006perspective}.
For more details on the implementation in SCIP, see~\cite{BestuzhevaGleixnerVigerske2021}.
An example is given in Figure~\ref{fig:persp}.

\subsection{Optimization-Based Bound Tightening}

Optimization-Based Bound Tightening (OBBT) is a domain propagation technique which minimizes and maximizes each variable over the feasible set of the problem or a relaxation thereof~\cite{QuGr93}.
Whereas FBBT (see Section~\ref{sec:consnl}) propagates the nonlinearities individually, OBBT considers (a relaxation of) all constraints together, and may hence compute tighter bounds.
However, it is rather expensive compared to FBBT.

In SCIP, OBBT solves two auxiliary LPs for each variable $x_k$ that could be subject to spatial branching:
\begin{equation}
  \label{eq:OBBTLP}
  \min / \max \{ x_k : D^xx+D^ww \leq d, c^\top x \leq U, x \in [\low{x},\upp{x}], w\in [\low{w},\upp{w}] \}
\end{equation}
where $D^xx+D^ww \leq d$, $D^x\in\mathbb{R}^{\ell \times n}$, $D^w\in\mathbb{R}^{\ell\times\hat m}$, $d\in\mathbb{R}^{\ell}$ is the linear relaxation of the feasible region of~\eqref{eq:minlp_ext}, and~$c^\top x \leq U$ is an objective cutoff constraint that excludes solutions with objective value worse than the current incumbent.
The optimal value of \eqref{eq:OBBTLP} may then be used to tighten the lower / upper bound of variable $x_k$.
A variable is subject to spatial branching if cut separation routines use the bounds of the variable at a node of the branch-and-bound tree.

SCIP, by default, applies OBBT at the root node to tighten bounds globally.
It restricts the computational effort by limiting the amount of LP iterations spent for solving the auxiliary LPs and interrupting for cheaper domain propagation techniques to be called between LP solves.

Further, SCIP does not only use the optimal objective values of~\eqref{eq:OBBTLP} to tighten the bounds on $x_k$, but it also applies a computationally cheap approximation of OBBT during the branch-and-bound search by exploiting the dual solutions from solves of~\eqref{eq:OBBTLP} at the root node.
Suppose the maximization LP is solved and feasible dual multipliers $\lambda_1,\ldots,\lambda_\ell$, $\mu \geq 0$ for $D^xx+D^ww \leq d$, $c^\top x \leq U$, respectively, and the corresponding reduced cost vectors $r^x$ and $r^w$ are obtained.
Then
\begin{equation}\label{eq:lvb}
  x_k \leq \sum_j r^x_j x_j + \sum_j r^w_j w_j + \lambda^\top d + \mu U
\end{equation}
is a valid inequality, which is called \emph{Lagrangian variable bound} (LVB), and
\begin{equation}\label{eq:lvbvalue}
  \sum_{j:r^x_j < 0} r^x_j \low{x}_j + \sum_{j:r^x_j > 0} r^x_j \upp{x}_j +
  \sum_{j:r^w_j < 0} r^w_j \low{w}_j + \sum_{j:r^w_j > 0} r^w_j \upp{w}_j + \lambda^\top d + \mu U
\end{equation}
is a valid upper bound for $x_k$ that equals the OBBT bound if the dual multipliers are optimal.
SCIP learns LVBs at the root node and propagates them during the tree search whenever the bounds of variables on the right-hand side of~\eqref{eq:lvb} become tighter or an improved primal solution is found.
For further details, see~\cite{GleixnerBertholdMuellerWeltge2017}.

In addition to OBBT with respect to the LP relaxation, also a variant is available that optimizes single variables over the potentially tighter convex NLP relaxation that is given by all linear and convex nonlinear constraints of~\eqref{eq:minlp}.
Also for this variant, linear Lagrangian variable bounds similar to~\eqref{eq:lvb} can be constructed by taking constraint convexity and KKT conditions into account.
Because of the potentially high computational cost of solving many NLPs, this variant of OBBT is deactivated by default.
For more details, see~\cite{SCIPoptsuite40}.

\subsection{Primal Heuristics}
\label{sec:heur}

The purpose of primal heuristics is to find high quality feasible solutions early in the search.
When given an MINLP, up to 40 primal heuristics are active in SCIP by default.
Many of them aim to find an integer-feasible solution to the LP relaxation.
In the following, primal heuristics that are only active in the presence of nonlinear constraints are discussed.

\subsubsection{subNLP}

A primal heuristic like \texttt{subNLP} is implemented in virtually any global MINLP solver.
Given a point $\tilde x$ that satisfies the integrality requirements ($\tilde x_{\mathcal{I}}\in\mathbb{Z}^{\vert\mathcal{I}\vert}$), the heuristic starts by fixing all integer variables in \eqref{eq:minlp} to the values given by $\tilde x$.
It then calls the SCIP presolver on this subproblem for possible simplifications.
Finally, it triggers a solution of the remaining NLP, using $\tilde x$ as the starting point.
If the NLP solver, such as Ipopt, finds a solution that is feasible (and often also locally optimal) for the NLP relaxation, then a feasible point for~\eqref{eq:minlp} has been found.

The starting point $\tilde x$ can be the current solution of the LP relaxation if integer-feasible, a point found by a primal heuristic that searches for integer-feasible solutions of the LP relaxation, or a point that is passed on by other primal heuristics for MINLP, such as those mentioned in the next sections.

How frequently the heuristic should run and how much effort to spend on an NLP solve is a nontrivial decision.
In the current implementation, the heuristic uses a fixed number for the iteration limit of the NLP solver for its first run.
For the following calls, the limit is set to twice the average number of iterations required in previous runs.
If, however, many of the previous runs hit the iteration limit, then an increased iteration limit is used.
Whether to run the heuristic at a node of the branch-and-bound tree depends on the number of nodes processed since it ran the last time, the iteration limit that would be used, and how successful the heuristic has been in finding feasible points in previous calls.

\subsubsection{Multistart}

If~\eqref{eq:minlp} is nonconvex after fixing all integer variables, then several local optima may exist for the NLPs solved by heuristic \texttt{subNLP}.
The success of the NLP solver then strongly depends on the starting point.
Therefore, the multistart heuristic aims to compute several starting points that are passed to the \texttt{subNLP} heuristic.

The algorithm, originally developed in~\cite{SmithChinneckAitken2013}, tries to approximate the boundary of the feasible set of the NLP relaxation by sampling points from $[\low{x},\upp{x}]$ and pushing them towards the feasible set by the use of an inexpensive gradient descent method.
Afterwards, points that are relatively close to each other are grouped into clusters.
Ideally, each cluster approximates the boundary of some connected component of the feasible set.
For each cluster, a linear combination of the points is passed as a starting point to \texttt{subNLP}.
For integer variables $x_i$, $i\in\mathcal{I}$, the value in the starting point is rounded to an integral value.

To reduce infeasibility of a point $\hat x$, the \emph{constraint consensus} method~\cite{SmithChinneckAitken2013} is used.
The algorithm computes a descent direction for each violated constraint of~\eqref{eq:minlp}.
For example, if $g_i(\hat x) > \upp{g}_i$ for some $i\in\{1,\ldots,m\}$, then the descent direction is given by $-\frac{g_i(\hat x)}{\Vert\nabla g_i(\hat x)\Vert^2} \nabla g_i(\hat x)$.
Point $\hat x$ is then updated by adding the average of the descent directions for all violated linear and nonlinear constraints.
This step is iterated until $\hat x$ becomes feasible, or a stopping criterion has been fulfilled.

The multistart heuristic currently runs for continuous problems ($\mathcal{I}=\emptyset$) only by default, since rounding and fixing integer variables most likely lead to infeasible NLP subproblems.
For more details, see~\cite{SCIPoptsuite40}.

\subsubsection{NLP Diving}

As an alternative to finding a good fixing for all integer variables of~\eqref{eq:minlp}, the NLP diving heuristic starts by solving the NLP relaxation at the current branch-and-bound node with an NLP solver.
It then iteratively fixes integer variables with fractional value and resolves both the LP and NLP relaxations, thereby simulating a depth-first-search in a branch-and-bound tree.
By default, variables for which the sum of the distances from the solutions of the LP and NLP relaxations to a common integer value is minimal are rounded to the nearest integer value.
Further, binary variables and nonlinear variables are preferred.
If the resulting NLP is found to be (locally) infeasible, one-level backtracking is applied, that is, the last fixing is undone, and the opposite fixing is tried.
If this is infeasible, too, the heuristic aborts.

\subsubsection{MPEC}

While the NLP diving heuristic either completely omits or enforces integrality restrictions in the NLP relaxation, the MPEC heuristic adds a relaxation of the integrality restriction to the NLP and tightens this relaxation iteratively.
The heuristic is only applicable to mixed-binary nonlinear programs at the moment.

The basic idea of the heuristic, originally developed in~\cite{ScheweSchmidt2019}, is to reformulate \eqref{eq:minlp} as a mathematical program with equilibrium constraints (MPEC) and to solve this MPEC to local optimality.
The MPEC is obtained from \eqref{eq:minlp} by rewriting the condition $x_i\in\{0,1\}$, $i\in\mathcal{I}$, as complementarity constraint $x_i \perp 1 - x_i$.
This reformulation is again reformulated to an NLP by writing it as $x_i\, (1-x_i) = 0$.
However, since these reformulated complementarity constraints will not, in general, satisfy constraint qualifications, solving this NLP reformulation with a generic NLP solver will often fail.

Therefore, in order to increase the chances of solving the NLP reformulation, the heuristic solves regularized versions of the NLP by relaxing $x_i (1-x_i) = 0$ to $x_i (1-x_i) \leq \theta$, for different, ever smaller $\theta > 0$.
The solution of one NLP is thereby used as the starting point for the next solve.
If the NLP solution is close to satisfying $x_\mathcal{I}\in\{0,1\}^{\vert\mathcal{I}\vert}$, it is passed as starting point to the \texttt{subNLP} heuristic.
If an NLP is (locally) infeasible, the heuristic does two more attempts where the values for binary variables that are already close to $0$ or $1$ are flipped to $1$ or $0$, respectively.
For more details, see~\cite{SCIPoptsuite50}.

\subsubsection{Undercover}

While the previous heuristics focused mainly on enforcing the integrality condition on an NLP, heuristic \texttt{undercover}~\cite{BertholdGleixner2014} starts from a completely different angle.
The heuristic is based on the observation that it sometimes suffices to fix only a comparatively small number of variables of~\eqref{eq:minlp} to yield a subproblem with all constraints being linear.
For example, for a bilinear term, only one of the variables needs to be fixed.
The variables to fix are chosen by solving a set covering problem, which aims at minimizing the number of variables to fix.
The values for the fixed variables are taken from the solution of the LP or NLP relaxation or a known feasible solution of the MINLP.

The resulting sub-MIP is less complex to solve, and does not need to be solved to proven optimality.
The solutions of the sub-MIP are immediately feasible for~\eqref{eq:minlp}.
However, the best one is also passed as starting point to heuristic \texttt{subnlp} to try for further improvement.
For more details, see~\cite{BertholdGleixner2014}.

\section{Benchmark}

The following aims to present a fair comparison of SCIP with several other state-of-the-art solvers for general MINLP.
Doing so is not trivial at all.
First, a set of instances needs to be selected that is suitable as a benchmark set.
Second, solver parameters have to be set such that all solvers solve the same instances with the same working limits and the same requirements on feasibility and optimality -- this goal could not be reached completely.
Third, the solver's results have to be checked for correctness, or, when this is not possible, plausibility.

GAMS was used for the experiments, as it provides various facilities to help on solver comparisons and comes with current versions of SCIP and the commercial solvers BARON~\cite{KhajaviradSahinidis2018}, Lindo~API~\cite{LinSchrage2009}, and Octeract included.
ANTIGONE has not been included in the comparison, as its development seems to have stopped years ago. %

All computations were run on a Linux cluster with Intel Xeon E5-2670~v2 CPUs (20 cores).
The GAMS version is 41.2.0, which includes SCIP~8.0.2, BARON~22.9.30, Lindo~API~14.0.5099.162, and Octeract~4.5.1.
A GAMS license with all solvers enabled was used, so that SCIP uses CPLEX~22.1.0.0 as LP solver and Ipopt with HSL MA27 as NLP solver, BARON can choose between all LP/MIP/NLP solvers that it interfaces with, and Octeract uses CPLEX~22.1.0.0 as LP/MIP/QP/QCP solver.

\subsection{Test Set}
\label{sec:testset}

To construct a test set suitable for benchmarking, the MINLPLib~\cite{minlplib} collection of 1595 MINLP instances was used as source.
First, all instances that could not be handled by some of the considered solvers were excluded, e.g., instances with trigonometric functions, as they are not supported by BARON.
All solvers were then run in serial mode (that is, with parallelization features disabled) on the remaining 1505 instances and using the parameter settings described below. %
The results of these runs were then used to select a set of 200 instances that could be solved by at least one solver, that were not all trivial, had a varying degree of integrality and nonlinearity, and such that having many instances with a similar name is avoided.
The latter was done to avoid overrepresentation of optimization problems for which many instances were added to MINLPLib.

Since small changes to an instance can lead to large variations in the solver's performance, the benchmark's reliability is improved by considering for each instance four additional variants where the order of variables and equations has been permuted.
The permuted instances were generated with GAMS/Convert.
Thus, a test set of 1000 instances is obtained.

The following approach was used to select the 200 instances before permutation:
Let $I$ be the set of 1505 instances, $d_i$ be the fraction of integer variables in instance $i\in I$, and $e_i$ be the fraction of nonzeros in the Jacobian and objective function gradient that correspond to nonlinear terms.
Next, assign to each instance an identifier $f_i\in F$ such that instances that seem to come from the same model are assigned the same identifier.
This goal is approximated by mapping~$i$ to the name of the instance until the first digit, underscore, or dash,
except for the block layout design instances \texttt{fo*}, \texttt{m*}, \texttt{no*}, \texttt{o*}, which were all assigned to the same identifier.
$\vert F\vert=230$ different identifiers were found this way.

Further, let $\upp{t}_i$ be the largest time in seconds that any solver who did not produce wrong results on instance $i$ spend on instance $i$.
Finally, let $S$ be the number of instances that could be solved by at least one solver.

To ensure that instances with a varying amount of integer variables and nonlinearity are included, the interval $[0,1]$ was split once at breakpoints $0.05, 0.25, 0.5, 0.9$ and once at $0.1, 0.25, 0.5$.
Let $D$ and $E$ be the resulting partitions of $[0,1]$.
For every interval from $D$ and $E$, the aim is to have roughly the same number of instances with $d_i$ and $e_i$ in the respective intervals.
For the choice of breakpoints that define $D$ and $E$, the distribution of $d_i$ and $e_i$, $i\in I$, have been taken into account.
For example, MINLPLib contains many purely continuous and purely discrete instances, but not many instances that are mostly linear or completely nonlinear.

To avoid including too many instances originating from the same model, including more than two instances for each identifier in $F$ is discouraged.
Further, instances that seem trivial, i.e., which are solved by all solvers in no more than five seconds, or could not be solved by any solver are excluded.
Introducing penalty terms, the following optimization problem for instance selection is obtained:
\begin{align*}
  \min\; & \sum_{d \in D} \lambda_d^2
         + \sum_{e \in E} \lambda_e^2
         + 10 \sum_{f\in F} \lambda_f^2 \\
 \text{such~that}\;
   & \sum_{i\in I: d_i\in d} z_i = \left\lfloor \frac{N}{\vert D\vert} \right\rceil + \lambda_{d} && \forall d \in D, \\
   & \sum_{i\in I: e_i \in e} z_i = \left\lfloor \frac{N}{\vert E\vert} \right\rceil + \lambda_{e} && \forall e \in E, \\
   & \sum_{i\in I: f_i = f} z_i \leq 2 + \lambda_f && \forall f \in F, \\
   & z_i = 0 && \forall i\in I: \upp{t}_i \leq 5, \\
   & z_i = 0 && \forall i\in I: i\not\in S, \\
   & z \in \{0,1\}^{\vert I\vert}, \lambda \in \mathbb{Z}^{\vert D\vert+\vert E\vert+\vert F\vert}
\end{align*}
This problem was solved for $N$ varying between 180 and 220.
For $N=208$, this yield a selection of 200 instances with an acceptable penalty value of 106.
See Section~\ref{sec:testsetdetail} for a list of all selected instances.
Table~\ref{tab:bucketsizes} shows the number of instances for each element of $D\times E$.
For five identifiers from $F$, three instead of two instances were selected, i.e., $\lambda_f=1$ for five $f\in F$.

\begin{table}[ht]
\centering
\begin{tabular}{l|rrrrr|r}
\toprule
$E\downarrow$ \,$\vert$\, $D\rightarrow$ & $[0,0.05)$ & $[0.05,0.25)$ & $[0.25,0.5)$ & $[0.5,0.9)$ & $[0.9,1]$  & $[0,1]$ \\ \midrule
$[0,0.1)$    & 3 & 7 & 19 & 15 & 6 & 50 \\
$[0.1,0.25)$ & 8 & 22 & 9 & 7 & 4 & 50 \\
$[0.25,0.5)$ & 8 & 8 & 6 & 10 & 18 & 50 \\
$[0.5,1]$    & 25 & 2 & 5 & 7 & 11 & 50 \\\midrule
$[0,1]$      & 44 & 39 & 39 & 39 & 39 & 200 \\
\bottomrule
\end{tabular}
\caption{Number of instances selected with ``discreteness'' $d_i$ and ``nonlinearity'' $e_i$ in intervals from $D$ and $E$.}
\label{tab:bucketsizes}
\end{table}

\subsection{Parameter Settings}

\subsubsection{Missing Variable Bounds}
\label{sec:missingbounds}

To compute a lower bound on the optimal value of a minimization problem, all solvers considered here construct a convex relaxation of the given problem.
For nonconvex constraints, this often relies on the computation of valid convex underestimators or concave overestimators.
As these typically depend on variables' bounds (recall the McCormick underestimators~\eqref{eq:mccormick}), missing or very large bounds on variables in nonconvex terms can mean that an instance will be very hard or impossible to solve.

Even when the user forgot to specify some variable bounds, the solver may still be able to derive bounds via domain propagation.
Further, once a feasible solution $\hat x$ has been found, additional bounds may be derived from the inequality $c^\top x\leq c^\top \hat x$.
However, as there are always cases where bounds are still missing after presolve, solvers invented different ways to deal with this obstacle.

If SCIP cannot construct an under- or overestimator because of missing variable bounds, it continues by branching on an unbounded variable.
This way, there will eventually be a node in the branch-and-bound tree where all variables are bounded.
Nodes that still contain unbounded variable domains may be pruned due to a derived lower bound on the objective function exceeding the incumbents objective function value.
But it may also be the case that pruning will not be possible and SCIP does not terminate.
However, variable bounds after branching cannot grow indefinitely in SCIP, but are limited by $\pm 10^{20}$ by default.
That is, SCIP does not search for solutions with variable values beyond this value.

The other solvers considered here add variable bounds based on a heuristic decision.
If BARON is still missing bounds on variables in nonconvex terms after presolve, it sets the bound to a value that depends on the type of nonlinearity involved.
Typically, this value is around $\pm 10^{10}$.
BARON also prints a warning to the log and no longer claim to have solved a problem to global optimality, i.e., it does not return a lower bound.
Lindo API adjusts the bounds for all variables that are involved in convexification to be within $[-10^{10},10^{10}]$.
At termination, it returns the lower bound for the restricted problem.
Octeract proceeds similarly and introduces a bound of $\pm 10^7$ for every missing bound and returns the lower bound for the restricted problem at termination.

Evidently, just passing an instance with unbounded variables to a solver with default settings may mean that each solver solves a different subproblem of the actual problem and often also reports a lower bound that corresponds to the solved subproblem only.
Fortunately, for every solver considered here, parameters are available to adjust the treatment of unbounded variables.
A first impulse could be to tell all solvers to set missing bounds to infinity, but this is not possible as each solver treats values beyond a certain finite value as ``infinity'' (BARON: $10^{50}$, Octeract: $10^{308}$, SCIP: $10^{20}$).
Changing this value is either not possible or not advisable.

We therefore decided to aim for $\pm 10^{12}$ as replacement for a missing variable bound.
For BARON and SCIP, the GAMS interface can replace any missing bound by $\pm 10^{12}$ before the instance is passed to the solver.
BARON will hence also return a lower bound for this restricted problem.
For Lindo API, a solver parameter can be changed so that bounds for all variables subject to convexification are bounded by $\pm 10^{12}$ (instead of $\pm10^{10}$).
Finally, also for Octeract, all missing bounds are set to $\pm10^{12}$ (instead of $\pm10^7$) by changing of a solver parameter.
Note, that this still does not ensure that all solvers solve the same instance, since Lindo API would still change initial finite bounds beyond $10^{12}$ and may also not set any bounds for variables that are not involved in convexification.

Next to missing bounds on problem variables, also singularities in functions (e.g., $1/x$, $\log(x)$) can prevent finite under- or overestimators from being available.
Unfortunately, there are no parameters available to ensure a uniform treatment of this case in all solvers.
SCIP ensures that a variable $x$ in $x^p$, $p<0$, or $\log(x)$ is bounded away from zero by $10^{-9}$, and terminates with a lower bound for this modified problem.
BARON applies the same method as the one for missing bounds on problem variables to choose a suitable bound on $x$.
No lower bound is returned at termination then.
The methods in Lindo API and Octeract are not known to us.

\subsubsection{Solution Quality}
\label{sub:tolerances}

To ensure that all solvers return solutions of the same quality, constraints of~\eqref{eq:minlp} are required to be satisfied with an absolute tolerance of $10^{-6}$. This applies to linear and nonlinear equations, variable bounds, and integrality.

In addition, a tolerance on the proof of optimality is set.
For this purpose, typically, solvers are allowed to stop when the absolute or relative gap between lower and upper bounds on the optimal value are sufficiently small.
Since the test set is diverse and has optimal values of varying magnitude, setting only a relative gap limit and no absolute gap limit would be preferable.
Unfortunately, Octeract does not permit different values for these limits.
As a compromise, BARON, Lindo API, and SCIP are run with $10^{-4}$ as relative gap limit and $10^{-6}$ as absolute gap limit, while for Octeract, $10^{-6}$ is used for both the absolute and relative gap limit.
Below, the impact of using a tighter optimality tolerance for Octeract is analyzed in a separate comparison.

\subsubsection{Working Limits}

As working limits, a time limit of two hours is used and the jobs on the cluster are restricted to 50 GB of RAM.
Further, the amount parallelization (multiple threads or processes) that a solver is allowed to use is limited in varying degrees.
To simplify the presentation, the term ``threads'' is used also for Octeract, even though it uses multiple processes instead of threads to parallelize its solving process.

\subsubsection{Summary}

To summarize, the following parameters are used:
\begin{description}
\item[GAMS] (applied to all solvers): \texttt{optcr=1e-4}, \texttt{optca=1e-6}, \texttt{reslim=7200}, \texttt{workspace={\allowbreak}50000}, \texttt{threads} $\in\{1,4,8,16\}$
\item[BARON:] \texttt{InfBnd=1e12}, \texttt{AbsConFeasTol=1e-6}, \texttt{AbsIntFeasTol=1e-6}
\item[Lindo API:] \texttt{GOP\_BNDLIM=1e12}, \texttt{SOLVER\_FEASTOL=1e-6}
\item[Octeract:] \texttt{INFINITY=1e12}, \texttt{INTEGRALITY\_VIOLATION\_TOLERANCE=1e-6}
\item[SCIP:] \texttt{gams/infbound=1e12}, \texttt{constraints/nonlinear/linearizeheursol=o} (this\\ undoes a change in the algorithmic settings of SCIP that is part of the GAMS/SCIP interface)
\end{description}

\subsection{Correctness Checks}

The GAMS/Examiner 2.0 tool is used to evaluate the violation of constraints, bounds, and integrality in the solutions reported by the solver.
Examiner generates for each solver a file that contains for each instance the solving time, returned lower and upper bound, and solution infeasibility.

A run of a solver on an instance is marked as \emph{failed} if the solver terminated abnormally, the solution is not feasible with respect to the feasibility tolerance, or the lower or upper bound contradicts with the bounds on the optimal value that are specified on the MINLPLib page.
Note, that the primal and dual bounds on the MINLPLib page were calculated without enforcing the $\pm 10^{12}$ limit on unbounded variables.
However, in order for an instance to be accepted into the test set, one of the solvers considered here must have solved the instance and found an optimal value that fits within the lower and upper bounds given at MINLPLib.
It is therefore acceptable to use these bounds for checking.

A run that has not failed is marked as \emph{solved} if the relative or absolute limits on the gap between lower and upper bound are satisfied.
If a solver stopped without closing the gap before the time limit, then the solver time is changed to the time limit.
The only exception here is BARON, which stops on two instances before the time limit without reporting a lower bound due to singularities in functions (see Section~\ref{sec:missingbounds}).
To be consistent with the treatment of other solvers, these two instances were accounted as solved by BARON with the original solver time.

\subsection{Results}

\subsubsection{Serial Mode}

For the main comparison, all parallelization features in the solvers were disabled, that is, GAMS was run with option \texttt{threads} set to 1.
In addition to the solver itself, results for the \emph{virtual best} and \emph{virtual worst} solver are reported, which are obtained by picking for each instance the fastest or the slowest solver, respectively.

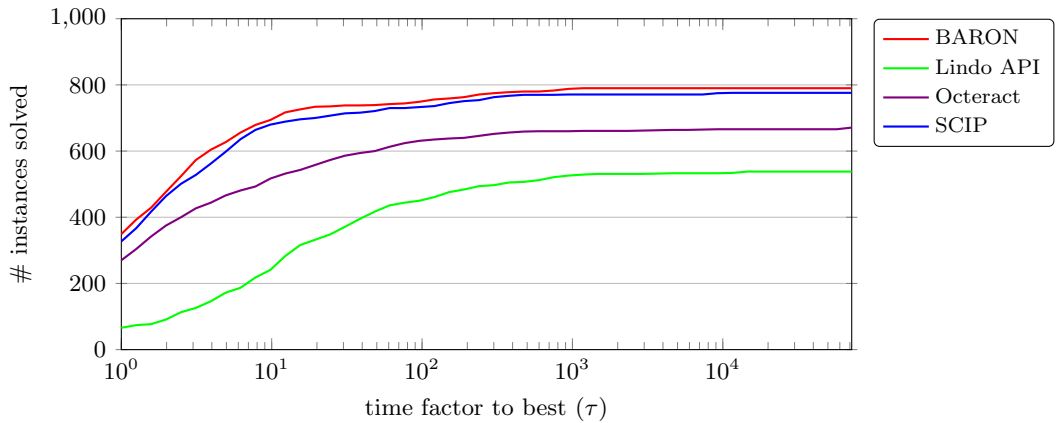
\begin{figure}[ht]
 \pgfplotsset{ legend style = { font = \footnotesize, rounded corners = 2pt } }
 \begin{tikzpicture}[font=\small]
  \begin{semilogxaxis}[
      width = .75\linewidth, height = 0.4\linewidth,
      xmin = 1, xmax = 72000,
      xlabel = { time factor to best ($\tau$)},
      axis y line = left,
      axis y line* = box,
      ylabel = {\# instances solved},
      ymin = 0, ymax = 1000,
      ytick = {0, 200, 400, 600, 800, 1000},
      ymajorgrids,
      legend pos = outer north east,
      legend cell align = left
      ]

      \addplot[red,thick] table[x=tau, y=BARON] {benchmark/pprofile.txt}
      ;%
      \addlegendentry{BARON};

      \addplot[green,thick] table[x=tau, y=Lindo] {benchmark/pprofile.txt}
      ;%
      \addlegendentry{Lindo API};

      \addplot[violet,thick] table[x=tau, y=Octeract] {benchmark/pprofile.txt}
      ;%
      \addlegendentry{Octeract};

      \addplot[blue,thick] table[x=tau, y=SCIP] {benchmark/pprofile.txt}
      ;%
      \addlegendentry{SCIP};

  \end{semilogxaxis}
 \end{tikzpicture}
 \caption{Performance profile comparing all solvers with parallelization disabled.}
 \label{fig:pprofile_1th}
\end{figure}

Table~\ref{tab:results_1th} shows for each solver the number of instances that could be solved, the number of times the time limit was reached, and the number of runs that were marked as failed.
Further, the shifted geometric mean of the running time of the solver is provided.
The shift has been set to 1 second.
Here, instances that failed are accounted with the time limit.
The performance profile~\cite{DolanMore2002} in Figure~\ref{fig:pprofile_1th} shows the number of instances a solver solved within a time that is at most a factor of the fastest solvers time.
Section~\ref{sec:detailed_singlethread} provides detailed results.

\begin{table}[ht]
 \centering
 \begin{tabular}{l|rrrr}
& solved & timeout & fail & time \\ \midrule
BARON & 790 & 183 & 27 & 75.4 \\
Lindo API & 538 & 323 & 139 & 489.1 \\
Octeract & 671 & 279 & 50 & 184.1 \\
SCIP & 776 & 183 & 41 & 85.2 \\
\midrule
virt.~worst & 368 & 405 & 227 & 1505.2 \\
virt.~best & 967 & 33 & 0 & 19.7 \\
\bottomrule
\end{tabular}

 \caption{Aggregated performance data for all solvers on test set of 1000 instances with parallelization disabled.}
 \label{tab:results_1th}
\end{table}

The results show a small lead of BARON before SCIP with respect to both number of instances solved and average time.
Since the number of timeouts is almost equal, one could argue that it is the higher stability of BARON that moves it onto the first place here.
In fact, the 41 fails of SCIP are due to returning a wrong optimal value 16 times, returning an infeasible solution 23 times, and aborts due to numerical troubles for two instances.
For BARON, fails are due to returning a wrong optimal value 26 times and an infeasible solution only once. %
While SCIP~8.0 has made a large step forward in ensuring that nonlinear constraints are satisfied in the non-presolved problem, violations in linear constraints or variable bounds still occur for a few instances.
These are typically due to variables being aggregated during presolve.

Even though Octeract and Lindo API solved considerably fewer instances than BARON and SCIP, which also results in an increased mean time, it is noteworthy that each of the two is also the fastest solver on 270 and 66 instances, respectively.
Octeract also produced correct results for 95\% of the test set, while for Lindo API a relatively high number of wrong optimal values, infeasible solutions, or aborts is observed.

The large differences between the real and virtual solvers show that none of the solvers dominates all others or is dominated.

Next, the effect of changing the gap limit for Octeract has been investigated.
Recall from Section~\ref{sub:tolerances} that relative and absolute gap limits of $10^{-6}$ and $10^{-4}$, respectively, were used for all solvers except for Octeract.
Since Octeract does not allow choosing these limits separately, it had been run with the tighter relative gap limit of $10^{-6}$.
To check whether this lead to a considerable disadvantage for this solver, the solver was rerun on the 200 non-permuted instances with both relative and absolute gap limit set to $10^{-4}$.
The table in Section~\ref{sec:detailed_octeractconvtol} shows that the change in the convergence tolerance had essentially no effect on the solver's performance.
In both cases, the same 134 instances could be solved.
The mean time changed from 178.6 for a limit of $10^{-6}$ to 179.0 for a limit of $10^{-4}$.

\subsubsection{Parallel Mode}

In the next comparison, each solver is allowed to use multiple threads or processes.
Since SCIP's use of multiple threads is limited to presolving MIPs, checking quadratic functions for convexity, and the linear algebra in Ipopt, FiberSCIP~\cite{ShinanoHeinzVigerskeWinkler2018} is used to run SCIP in parallel mode.
FiberSCIP is a shared-memory instantiation of the UG framework~\cite{ug} for the parallelization of branch-and-bound based solvers.
The framework parallelizes the search of the branch-and-bound tree by collecting and distributing open problems between independent instances of SCIP.
In addition, the first seconds of the solving process are used for a ``racing ramp-up'' phase.
Here, multiple SCIP instances with differing parameter sets are run concurrently, and the one with the best lower bound is used for the remaining solve.
The UG version was 1.0.0~beta3.

For the runs in serial mode, reaching the memory limit of 50 GB was not observed for any solver.
But since parallelization often increases memory requirements, a memory limit of 100 GB has been used for the runs in parallel mode.
Since this meant a reduction in available computing resources, only the 200 non-permuted instances are used for comparisons. %

\begin{figure}[ht]
 \pgfplotsset{ legend style = { font = \footnotesize, rounded corners = 2pt } }
 \begin{tikzpicture}[font=\small]
  \begin{semilogxaxis}[
      width = .75\linewidth, height = 0.4\linewidth,
      xmin = 1, xmax = 1800,
      xlabel = { time factor to best ($\tau$)},
      axis y line = left,
      axis y line* = box,
      ylabel = {\# instances solved},
      ymin = 0, ymax = 200,
      ytick = {0, 50, 100, 150, 200},
      ymajorgrids,
      legend pos = outer north east,
      legend cell align = left
      ]

      \addplot[red,thick] table[x=tau, y=SCIP.1] {benchmark/pprofile_scipmultithread.txt};
      \addlegendentry{1 thread};

      \addplot[green,thick] table[x=tau, y=SCIP.4] {benchmark/pprofile_scipmultithread.txt};
      \addlegendentry{4 threads};

      \addplot[violet,thick] table[x=tau, y=SCIP.8] {benchmark/pprofile_scipmultithread.txt};
      \addlegendentry{8 threads};

      \addplot[blue,thick] table[x=tau, y=SCIP.16] {benchmark/pprofile_scipmultithread.txt};
      \addlegendentry{16 threads};

      \addplot[black,dashed] table[x=tau, y=virtbest] {benchmark/pprofile_scipmultithread.txt};
      \addlegendentry{virt.\ best};

      \addplot[black,dashed] table[x=tau, y=virtworst] {benchmark/pprofile_scipmultithread.txt};
      \addlegendentry{virt.\ worst};

  \end{semilogxaxis}
 \end{tikzpicture}
 \caption{Performance profile comparing SCIP and FiberSCIP.}
 \label{fig:pprofile_fiberscip}
\end{figure}
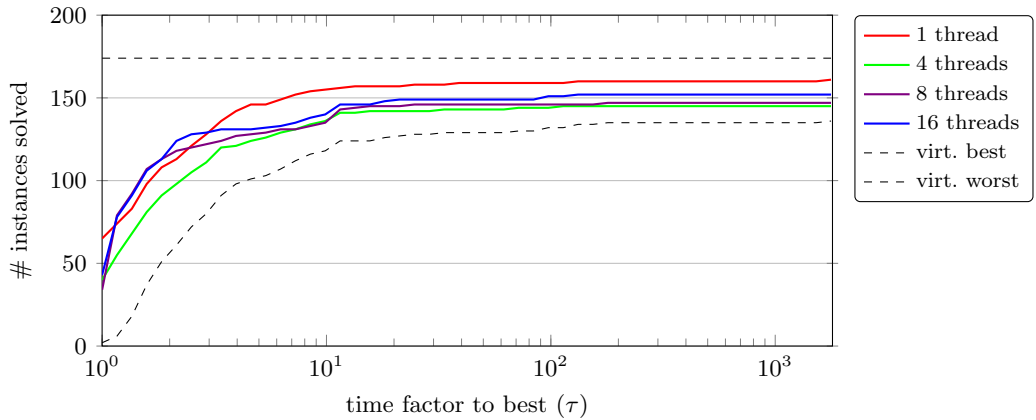

Table~\ref{tab:results_multithread} shows, for an increasing number of threads, the number of instances that could be solved by each solver and the mean time spent.
In addition, Figure~\ref{fig:pprofile_fiberscip} provides a performance profile that compares SCIP and FiberSCIP only.
Section~\ref{sec:detailed_multithread} gives detailed results.

\begin{table}[ht]
 \centering
 \begin{tabular}{l|rrrrrrrr}
& \multicolumn{2}{c}{1 thread} & \multicolumn{2}{c}{4 threads} & \multicolumn{2}{c}{8 threads} & \multicolumn{2}{c}{16 threads} \\
 & solved & time & solved & time & solved & time & solved & time \\ \midrule
BARON & 161 & 64.3 & 160 & 58.2 & 160 & 57.1 & 158 & 58.6 \\
Lindo API & 114 & 423.6 & 114 & 379.2 & 106 & 459.5 & 107 & 456.4 \\
Octeract & 134 & 178.6 & 133 & 146.9 & 138 & 118.1 & 135 & 123.2 \\
(Fiber)SCIP & 161 & 76.9 & 145 & 94.3 & 147 & 77.8 & 152 & 74.8 \\
\bottomrule
\end{tabular}

 \caption{Aggregated performance data for all solvers on test set of 200 instances when run with parallelization allowed.}
 \label{tab:results_multithread}
\end{table}

Apparently, enabling parallelization seldom has a considerable advantage on this test set.
For Octeract, where parallelization was part of its original design, a small increase in the number of instances that could be solved and a reduction in time by~34\% when using up to 8 parallel processes is observed.
As far as we know, BARON's use of multiple threads is currently limited to enabling this feature in the solver for a MIP relaxation.
As a consequence, only moderate improvement of running time by up to~11\% are seen.
For Lindo API, an improvement due to parallelization seems to be impeded by a further increase in fails when using multiple threads (1~thread: 24, 4: 28, 8: 35, 16: 43).
Finally, for SCIP/FiberSCIP the additional overhead due to the parallelization being build on top of the solver instead of being tightly integrated is not compensated by the use of multiple threads.
However, in contrast to other solvers, a monotonous improvement in both number of instances solved and mean solving time when increasing from 4 to 16 threads is observed.
Further, the virtual solvers in the performance profile show that FiberSCIP can solve instances that SCIP on one thread couldn't solve. %

Finally, note that a benefit due to parallelization can usually only be expected for rather challenging instances because of the additional overhead in duplicating and synchronizing data and processes.
However, the test set deliberately included only instances that could already be solved by some solver in serial mode, and only instances that were trivial for all solvers, though they may be solved quickly by some, were excluded.
As a small experiment, for each solver only those instances that required at least 10 or 100 seconds to solve in serial mode were considered.
Unfortunately, this essentially repeated the trends shown in Table~\ref{tab:results_multithread}, so details are omitted here. %
A more thorough analysis of the parallelization capabilities of MINLP solvers using a set of challenging instances only would be necessary, but exceeds the scope of this paper.

\section{Conclusion}

The development of the MINLP solver in SCIP has come a long way.
In a recent version-to-version comparison~\cite[slides 49-51]{marcscipversions}, a steady improvement in the performance of SCIP on MINLP over the last ten years has been measured, resulting in SCIP 8 solving twice as many instances as SCIP 3 and a speed-up of factor three.
Partially, this improvement has been achieved by improving and adding features particular for MINLP.
However, due to the generality of SCIP as a CIP solver, also many developments that targeted MIP solving were immediately available for MINLP solving.

With version 8.0, the MINLP solving capabilities of SCIP have been largely reworked and extended, which resulted in a considerable improvement in both robustness and performance~\cite{SCIPoptsuite80,marcscipversions}.
As a result, SCIP's performance is currently on par with the state-of-the-art commercial solver BARON.

In contrast to the commercial solvers considered here, SCIP offers a variety of possibilities for a user, developer, or researcher to interact with the solving process.
In particular, the newly added ``nonlinear handler'' plugin type sets SCIP apart from most other MINLP solvers, as it allows focusing on experimenting with new algorithms to handle certain structures in nonlinear functions without modifications to the solver's code.

The rather large number of features that are disabled by default shows that tuning and improving the existing code base has become increasingly necessary.
Future work will of course also include the addition of new features, e.g., improved separation for signomial functions~\cite{XuDAmbrosioLibertiVanier2022}, the use of alternative relaxations for polynomial functions~\cite{BestuzhevaGleixnerVoelker2022}, or monoidal strengthening of intersection cuts for quadratic constraints~\cite{ChmielaMunozSerrano2022}.

The increasing number of cores in present-day CPUs means that to fully utilize an ordinary desktop computer, a solver needs to be parallelized.
While the UG framework provides such a possibility for SCIP in both shared and distributed memory environments, the experiments with FiberSCIP on up to 16 threads show that more tuning is necessary to ensure that the additional overhead can be compensated by the use of additional computing resources.
Since the development of UG was initially motivated and has focused primarily on the use of large-scale parallel computing environments~\cite{ShinanoAchterbergBertholdHeinzKochWinkler2016}, an investigation on using UG with SCIP to solve challenging MINLPs in distributed memory environments with many CPU cores could be interesting as well.

\section*{Acknowledgments}

We are very much in all SCIP developers' debt -- the extensions to support nonlinear constraints and solve MINLPs would not have been possible without the framework's existence and the powerful MIP solver that we could build upon.
While the authors of this paper are the main developers of the new MINLP features in SCIP 8, many have contributed to the MINLP capabilities in previous releases of SCIP, namely Martin Ballerstein, Timo Berthold, Tobias Fischer, Thorsten Gellermann, Ambros Gleixner, Renke Kuhlmann, Dennis Michaels, Marc Pfetsch, and Stefan Weltge.
Further, we thank Yuji Shinano for the development of FiberSCIP and swiftly reacting to our request for the possibility to set gap limits.
Last but not least, we are very grateful to Franziska Schl{\"o}sser for the setup and maintenance of benchmarking and testing facilities for the infamous ``consexpr'' development branch of SCIP.

The work for this article has been conducted within the Research Campus Modal funded by the German Federal Ministry of Education and Research (BMBF grant numbers 05M14ZAM, 05M20ZBM).
Additional funding has been received from the German Federal Ministry for Economic Affairs and Energy within the project EnBA-M (ID: 03ET1549D).

\bibliographystyle{abbrvnat}
\bibliography{consexpr}

\begin{thebibliography}{78}
\providecommand{\natexlab}[1]{#1}
\providecommand{\url}[1]{\texttt{#1}}
\expandafter\ifx\csname urlstyle\endcsname\relax
  \providecommand{\doi}[1]{doi: #1}\else
  \providecommand{\doi}{doi: \begingroup \urlstyle{rm}\Url}\fi

\bibitem[Achterberg(2007)]{Achterberg2007}
T.~Achterberg.
\newblock \emph{Constraint Integer Programming}.
\newblock PhD thesis, Technische Universit{\"a}t Berlin, 2007.

\bibitem[Adams and Sherali(1986)]{adams1986tight}
W.~P. Adams and H.~D. Sherali.
\newblock A tight linearization and an algorithm for zero-one quadratic
  programming problems.
\newblock \emph{Management Science}, 32\penalty0 (10):\penalty0 1274--1290,
  1986.
\newblock \doi{10.1287/mnsc.32.10.1274}.

\bibitem[Adams and Sherali(1990)]{adams1990linearization}
W.~P. Adams and H.~D. Sherali.
\newblock Linearization strategies for a class of zero-one mixed integer
  programming problems.
\newblock \emph{Operations Research}, 38\penalty0 (2):\penalty0 217--226, 1990.
\newblock \doi{10.1287/opre.38.2.217}.

\bibitem[Adjiman and Floudas(1996)]{AdjimanFloudas1996}
C.~Adjiman and C.~Floudas.
\newblock Rigorous convex underestimators for general twice-differentiable
  problems.
\newblock \emph{Journal of Global Optimization}, 9\penalty0 (1):\penalty0
  23--40, 1996.
\newblock \doi{10.1007/BF00121749}.

\bibitem[Balas(1971)]{Balas1971}
E.~Balas.
\newblock Intersection cuts -- a new type of cutting planes for integer
  programming.
\newblock \emph{Operations Research}, 19\penalty0 (1):\penalty0 19--39, 1971.
\newblock \doi{10.1287/opre.19.1.19}.

\bibitem[Ballerstein et~al.(2015)Ballerstein, Michaels, and Vigerske]{BaMiVi12}
M.~Ballerstein, D.~Michaels, and S.~Vigerske.
\newblock Linear underestimators for bivariate functions with a fixed convexity
  behavior.
\newblock ZIB-Report 13-02, Zuse Institute Berlin, 2015.
\newblock \urnlink{nbn:de:0297-zib-17641}.

\bibitem[Beale(1980)]{Beale1980}
E.~Beale.
\newblock Branch and bound methods for numerical optimization of non-convex
  functions.
\newblock In M.~Barritt and D.~Wishart, editors, \emph{COMPSTAT 80 Proceedings
  in Computational Statistics}, pages 11--20, Vienna, 1980. Physica-Verlag.

\bibitem[Bell()]{cppad}
B.~Bell.
\newblock {CppAD}: A package for differentiation of {C++} algorithms.
\newblock \url{https://github.com/coin-or/CppAD/}.

\bibitem[Belotti et~al.(2013)Belotti, Kirches, Leyffer, Linderoth, Luedtke, and
  Mahajan]{BeKiLeLiLuMa12}
P.~Belotti, C.~Kirches, S.~Leyffer, J.~Linderoth, J.~Luedtke, and A.~Mahajan.
\newblock Mixed-integer nonlinear optimization.
\newblock \emph{Acta Numerica}, 22:\penalty0 1--131, 2013.
\newblock \doi{10.1017/S0962492913000032}.

\bibitem[Berthold and Gleixner(2014)]{BertholdGleixner2014}
T.~Berthold and A.~Gleixner.
\newblock Undercover: a primal {MINLP} heuristic exploring a largest {sub-MIP}.
\newblock \emph{Mathematical Programming}, 144\penalty0 (1-2):\penalty0
  315--346, 2014.
\newblock \doi{10.1007/s10107-013-0635-2}.

\bibitem[Berthold et~al.(2009)Berthold, Heinz, and Pfetsch]{BerHP09}
T.~Berthold, S.~Heinz, and M.~E. Pfetsch.
\newblock Nonlinear pseudo-boolean optimization: relaxation or propagation?
\newblock In O.~Kullmann, editor, \emph{Theory and Applications of
  Satisfiability Testing -- SAT 2009}, number 5584 in LNCS, pages 441--446,
  Berlin, Heidelberg, 2009. Springer.
\newblock \doi{10.1007/978-3-642-02777-2\_40}.

\bibitem[Berthold et~al.(2011)Berthold, Heinz, Pfetsch, and
  Vigerske]{BeHePfVi11}
T.~Berthold, S.~Heinz, M.~E. Pfetsch, and S.~Vigerske.
\newblock Large neighborhood search beyond {MIP}.
\newblock In L.~D. Gaspero, A.~Schaerf, and T.~Stützle, editors,
  \emph{Proceedings of the 9th Metaheuristics International Conference (MIC
  2011)}, pages 51--60, 2011.

\bibitem[Berthold et~al.(2012)Berthold, Heinz, and Vigerske]{BeHeVi09}
T.~Berthold, S.~Heinz, and S.~Vigerske.
\newblock Extending a {CIP} framework to solve {MIQCPs}.
\newblock In J.~Lee and S.~Leyffer, editors, \emph{Mixed Integer Nonlinear
  Programming}, volume 154 of \emph{The IMA Volumes in Mathematics and its
  Applications}, pages 427--444. Springer New York, NY, 2012.
\newblock \doi{10.1007/978-1-4614-1927-3\_15}.

\bibitem[Bestuzheva et~al.(2021{\natexlab{a}})Bestuzheva, Besan{\c{c}}on, Chen,
  Chmiela, Donkiewicz, van Doornmalen, Eifler, Gaul, Gamrath, Gleixner,
  Gottwald, Graczyk, Halbig, Hoen, Hojny, van~der Hulst, Koch, L{\"u}bbecke,
  Maher, Matter, M{\"u}hmer, M{\"u}ller, Pfetsch, Rehfeldt, Schlein,
  Schl{\"o}sser, Serrano, Shinano, Sofranac, Turner, Vigerske, Wegscheider,
  Wellner, Weninger, and Witzig]{SCIPoptsuite80}
K.~Bestuzheva, M.~Besan{\c{c}}on, W.-K. Chen, A.~Chmiela, T.~Donkiewicz, J.~van
  Doornmalen, L.~Eifler, O.~Gaul, G.~Gamrath, A.~Gleixner, L.~Gottwald,
  C.~Graczyk, K.~Halbig, A.~Hoen, C.~Hojny, R.~van~der Hulst, T.~Koch,
  M.~L{\"u}bbecke, S.~J. Maher, F.~Matter, E.~M{\"u}hmer, B.~M{\"u}ller, M.~E.
  Pfetsch, D.~Rehfeldt, S.~Schlein, F.~Schl{\"o}sser, F.~Serrano, Y.~Shinano,
  B.~Sofranac, M.~Turner, S.~Vigerske, F.~Wegscheider, P.~Wellner, D.~Weninger,
  and J.~Witzig.
\newblock {The SCIP Optimization Suite 8.0}.
\newblock ZIB Report 21-41, Zuse Institute Berlin, 2021{\natexlab{a}}.
\newblock \urnlink{nbn:de:0297-zib-85309}.

\bibitem[Bestuzheva et~al.(2021{\natexlab{b}})Bestuzheva, Gleixner, and
  Vigerske]{BestuzhevaGleixnerVigerske2021}
K.~Bestuzheva, A.~Gleixner, and S.~Vigerske.
\newblock A computational study of perspective cuts.
\newblock ZIB Report 21-07, Zuse Institute Berlin, 2021{\natexlab{b}}.

\bibitem[Bestuzheva et~al.(2022{\natexlab{a}})Bestuzheva, Gleixner, and
  Achterberg]{AchterbergBestuzhevaGleixner2022}
K.~Bestuzheva, A.~Gleixner, and T.~Achterberg.
\newblock Efficient separation of {RLT} cuts for implicit and explicit bilinear
  products.
\newblock Technical Report 2211.13545, arXiv, 2022{\natexlab{a}}.
\newblock \doi{10.48550/ARXIV.2211.13545}.

\bibitem[Bestuzheva et~al.(2022{\natexlab{b}})Bestuzheva, Gleixner, and
  V\"olker]{BestuzhevaGleixnerVoelker2022}
K.~Bestuzheva, A.~Gleixner, and H.~V\"olker.
\newblock Strengthening {SONC} relaxations with constraints derived from
  variable bounds.
\newblock In preparation, 2022{\natexlab{b}}.

\bibitem[Bley et~al.(2008)Bley, Koch, and Niu]{BleyKochNiu2008}
A.~Bley, T.~Koch, and L.~Niu.
\newblock Experiments with nonlinear extensions to {SCIP}.
\newblock ZIB-Report 08-28, Zuse Institute Berlin, 2008.
\newblock \urnlink{nbn:de:0297-zib-8300}.

\bibitem[Bussieck and Vigerske(2010)]{BussieckVigerske2010}
M.~R. Bussieck and S.~Vigerske.
\newblock {MINLP} solver software.
\newblock In J.~J. Cochran, L.~A. Cox, Jr., P.~Keskinocak, J.~P. Kharoufeh, and
  J.~C. Smith, editors, \emph{Wiley Encyclopedia of Operations Research and
  Management Science}. Wiley \& Sons, Inc., NJ, 2010.
\newblock \doi{10.1002/9780470400531.eorms0527}.

\bibitem[Büskens and Wassel(2013)]{BueskensWassel2013}
C.~Büskens and D.~Wassel.
\newblock The {ESA} {NLP} solver {WORHP}.
\newblock In G.~Fasano and J.~D. Pintér, editors, \emph{Modeling and
  Optimization in Space Engineering}, volume~73 of \emph{Springer Optimization
  and Its Applications}, pages 85--110. Springer New York, NY, 2013.
\newblock \doi{10.1007/978-1-4614-4469-5\_4}.

\bibitem[Chmiela et~al.(2021)Chmiela, Mu{\~{n}}oz, and
  Serrano]{ChmielaMunozSerrano2021}
A.~Chmiela, G.~Mu{\~{n}}oz, and F.~Serrano.
\newblock On the implementation and strengthening of intersection cuts for
  {QCQPs}.
\newblock In M.~Singh and D.~P. Williamson, editors, \emph{Integer Programming
  and Combinatorial Optimization}, pages 134--147, Cham, 2021. Springer.
\newblock \doi{10.1007/978-3-030-73879-2\_10}.

\bibitem[Chmiela et~al.(2022)Chmiela, Mu\~noz, and
  Serrano]{ChmielaMunozSerrano2022}
A.~Chmiela, G.~Mu\~noz, and F.~Serrano.
\newblock Monoidal strengthening and unique lifting in {MIQCP}s.
\newblock Submitted, preprint available at
  \url{https://www.gonzalomunoz.org/publications/}, 2022.

\bibitem[Dolan and Mor{\'e}(2002)]{DolanMore2002}
E.~D. Dolan and J.~J. Mor{\'e}.
\newblock Benchmarking optimization software with performance profiles.
\newblock \emph{Mathematical Programming}, 91\penalty0 (2):\penalty0 201--213,
  2002.
\newblock \doi{10.1007/s101070100263}.

\bibitem[Domes and Neumaier(2010)]{DomesNeumaier2010}
F.~Domes and A.~Neumaier.
\newblock Constraint propagation on quadratic constraints.
\newblock \emph{Constraints}, 15\penalty0 (3):\penalty0 404--429, 2010.
\newblock \doi{10.1007/s10601-009-9076-1}.

\bibitem[Duran and Grossmann(1986)]{DuranGrossmann1986}
M.~A. Duran and I.~E. Grossmann.
\newblock An outer-approximation algorithm for a class of mixed-integer
  nonlinear programs.
\newblock \emph{Mathematical Programming}, 36\penalty0 (3):\penalty0 307--339,
  1986.
\newblock \doi{10.1007/BF02592064}.

\bibitem[Fischer(2017)]{Fischer2017}
T.~Fischer.
\newblock \emph{Branch-and-cut for complementarity and cardinality constrained
  linear programs}.
\newblock PhD thesis, Technical University of Darmstadt, 2017.

\bibitem[Fletcher and Leyffer(1998)]{FletcherLeyffer1998}
R.~Fletcher and S.~Leyffer.
\newblock User manual for {filterSQP}.
\newblock Numerical Analysis Report NA/181, Department of Mathematics,
  University of Dundee, Scotland, 1998.

\bibitem[Floudas(1995)]{Floudas1995}
C.~A. Floudas.
\newblock \emph{Nonlinear and Mixed Integer Optimization: Fundamentals and
  Applications}.
\newblock Oxford University Press, New York, 1995.

\bibitem[Frangioni and Gentile(2006)]{frangioni2006perspective}
A.~Frangioni and C.~Gentile.
\newblock Perspective cuts for a class of convex 0--1 mixed integer programs.
\newblock \emph{Mathematical Programming}, 106\penalty0 (2):\penalty0 225--236,
  2006.
\newblock \doi{10.1007/s10107-005-0594-3}.

\bibitem[Gamrath et~al.(2016)Gamrath, Fischer, Gally, Gleixner, Hendel, Koch,
  Maher, Miltenberger, M{\"u}ller, Pfetsch, Puchert, Rehfeldt, Schenker,
  Schwarz, Serrano, Shinano, Vigerske, Weninger, Winkler, Witt, and
  Witzig]{SCIPoptsuite32}
G.~Gamrath, T.~Fischer, T.~Gally, A.~Gleixner, G.~Hendel, T.~Koch, S.~J. Maher,
  M.~Miltenberger, B.~M{\"u}ller, M.~E. Pfetsch, C.~Puchert, D.~Rehfeldt,
  S.~Schenker, R.~Schwarz, F.~Serrano, Y.~Shinano, S.~Vigerske, D.~Weninger,
  M.~Winkler, J.~T. Witt, and J.~Witzig.
\newblock {The SCIP Optimization Suite 3.2}.
\newblock ZIB Report 15-60, Zuse Institute Berlin, 2016.
\newblock \urnlink{nbn:de:0297-zib-57675}.

\bibitem[Gamrath et~al.(2020)Gamrath, Anderson, Bestuzheva, Chen, Eifler,
  Gasse, Gemander, Gleixner, Gottwald, Halbig, Hendel, Hojny, Koch, Bodic,
  Maher, Matter, Miltenberger, M\"uhmer, M\"uller, Pfetsch, Schl\"osser,
  Serrano, Shinano, Tawfik, Vigerske, Wegscheider, Weninger, and
  Witzig]{SCIPoptsuite70}
G.~Gamrath, D.~Anderson, K.~Bestuzheva, W.-K. Chen, L.~Eifler, M.~Gasse,
  P.~Gemander, A.~Gleixner, L.~Gottwald, K.~Halbig, G.~Hendel, C.~Hojny,
  T.~Koch, P.~L. Bodic, S.~J. Maher, F.~Matter, M.~Miltenberger, E.~M\"uhmer,
  B.~M\"uller, M.~E. Pfetsch, F.~Schl\"osser, F.~Serrano, Y.~Shinano,
  C.~Tawfik, S.~Vigerske, F.~Wegscheider, D.~Weninger, and J.~Witzig.
\newblock {The SCIP Optimization Suite 7.0}.
\newblock ZIB Report 20-10, Zuse Institute Berlin, 2020.
\newblock \urnlink{nbn:de:0297-zib-78023}.

\bibitem[Gleixner et~al.(2012)Gleixner, Held, Huang, and Vigerske]{GlHeHuVi12}
A.~Gleixner, H.~Held, W.~Huang, and S.~Vigerske.
\newblock Towards globally optimal operation of water supply networks.
\newblock \emph{Numerical Algebra, Control and Optimization}, 2\penalty0
  (4):\penalty0 695--711, 2012.
\newblock \doi{10.3934/naco.2012.2.695}.

\bibitem[Gleixner et~al.(2017{\natexlab{a}})Gleixner, Berthold, M{\"u}ller, and
  Weltge]{GleixnerBertholdMuellerWeltge2017}
A.~Gleixner, T.~Berthold, B.~M{\"u}ller, and S.~Weltge.
\newblock Three enhancements for optimization-based bound tightening.
\newblock \emph{Journal of Global Optimization}, 67\penalty0 (4):\penalty0
  731--757, 2017{\natexlab{a}}.
\newblock \doi{10.1007/s10898-016-0450-4}.

\bibitem[Gleixner et~al.(2017{\natexlab{b}})Gleixner, Eifler, Gally, Gamrath,
  Gemander, Gottwald, Hendel, Hojny, Koch, Miltenberger, M{\"u}ller, Pfetsch,
  Puchert, Rehfeldt, Schl{\"o}sser, Serrano, Shinano, Viernickel, Vigerske,
  Weninger, Witt, and Witzig]{SCIPoptsuite50}
A.~Gleixner, L.~Eifler, T.~Gally, G.~Gamrath, P.~Gemander, R.~L. Gottwald,
  G.~Hendel, C.~Hojny, T.~Koch, M.~Miltenberger, B.~M{\"u}ller, M.~E. Pfetsch,
  C.~Puchert, D.~Rehfeldt, F.~Schl{\"o}sser, F.~Serrano, Y.~Shinano, J.~M.
  Viernickel, S.~Vigerske, D.~Weninger, J.~T. Witt, and J.~Witzig.
\newblock {The SCIP Optimization Suite 5.0}.
\newblock ZIB Report 17-61, Zuse Institute Berlin, 2017{\natexlab{b}}.
\newblock \urnlink{nbn:de:0297-zib-66297}.

\bibitem[Gleixner et~al.(2020)Gleixner, Maher, M{\"u}ller, and
  Pedroso]{GleixnerMaherMuellerPedroso2020}
A.~Gleixner, S.~J. Maher, B.~M{\"u}ller, and J.~P. Pedroso.
\newblock Price-and-verify: a new algorithm for recursive circle packing using
  {Dantzig--Wolfe} decomposition.
\newblock \emph{Annals of Operations Research}, 284\penalty0 (2):\penalty0
  527--555, 2020.
\newblock \doi{10.1007/s10479-018-3115-5}.

\bibitem[Glover(1974)]{Glover1974}
F.~Glover.
\newblock Polyhedral convexity cuts and negative edge extensions.
\newblock \emph{Zeitschrift für Operations Research}, 18:\penalty0 181--186,
  1974.
\newblock \doi{10.1007/BF02026599}.

\bibitem[Grossmann and Kravanja(1997)]{GrossmannKravanja1997}
I.~E. Grossmann and Z.~Kravanja.
\newblock Mixed-integer nonlinear programming: A survey of algorithms and
  applications.
\newblock In A.~R. Conn, L.~T. Biegler, T.~F. Coleman, and F.~N. Santosa,
  editors, \emph{Large-Scale Optimization with Applications, Part {II}: Optimal
  Design and Control}, volume~93 of \emph{The IMA Volumes in Mathematics and
  its Applications}, pages 73--100. Springer, New York, 1997.
\newblock \doi{10.1007/978-1-4612-1960-6\_5}.

\bibitem[Hansen et~al.(1993)Hansen, Jaumard, Ruiz, and
  Xiong]{HansenJaumardRuizXiong1993}
P.~Hansen, B.~Jaumard, M.~Ruiz, and J.~Xiong.
\newblock Global minimization of indefinite quadratic functions subject to box
  constraints.
\newblock \emph{Naval Research Logistics (NRL)}, 40\penalty0 (3):\penalty0
  373--392, 1993.
\newblock
  \doi{10.1002/1520-6750(199304)40:3<373::AID-NAV3220400307>3.0.CO;2-A}.

\bibitem[Hojny and Pfetsch(2019)]{HojnyPfetsch2019}
C.~Hojny and M.~E. Pfetsch.
\newblock Polytopes associated with symmetry handling.
\newblock \emph{Mathematical Programming}, 175\penalty0 (1):\penalty0 197--240,
  2019.
\newblock \doi{10.1007/s10107-018-1239-7}.

\bibitem[Khajavirad and Sahinidis(2018)]{KhajaviradSahinidis2018}
A.~Khajavirad and N.~V. Sahinidis.
\newblock A hybrid {LP/NLP} paradigm for global optimization relaxations.
\newblock \emph{Mathematical Programming Computation}, 10\penalty0
  (3):\penalty0 383--421, 2018.
\newblock \doi{10.1007/s12532-018-0138-5}.

\bibitem[Kocis and Grossmann(1989)]{KocisGrossmann1989}
G.~Kocis and I.~Grossmann.
\newblock Computational experience with {DICOPT}: Solving {MINLP} problems in
  process systems engineering.
\newblock \emph{Computers \& Chemical Engineering}, 13\penalty0 (3):\penalty0
  307--315, 1989.
\newblock \doi{10.1016/0098-1354(89)85008-2}.

\bibitem[Kronqvist et~al.(2016)Kronqvist, Lundell, and
  Westerlund]{KronqvistLundellWesterlund2016}
J.~Kronqvist, A.~Lundell, and T.~Westerlund.
\newblock The extended supporting hyperplane algorithm for convex mixed-integer
  nonlinear programming.
\newblock \emph{Journal of Global Optimization}, 64\penalty0 (2):\penalty0
  249--272, 2016.
\newblock \doi{10.1007/s10898-015-0322-3}.

\bibitem[Liberti(2012)]{Liberti2012a}
L.~Liberti.
\newblock Reformulations in mathematical programming: automatic symmetry
  detection and exploitation.
\newblock \emph{Mathematical Programming}, 131\penalty0 (1):\penalty0 273--304,
  2012.
\newblock \doi{10.1007/s10107-010-0351-0}.

\bibitem[Lin and Schrage(2009)]{LinSchrage2009}
Y.~Lin and L.~Schrage.
\newblock The global solver in the {LINDO API}.
\newblock \emph{Optimization Methods \& Software}, 24\penalty0 (4--5):\penalty0
  657--668, 2009.
\newblock \doi{10.1080/10556780902753221}.

\bibitem[Locatelli(2018)]{Locatelli2018}
M.~Locatelli.
\newblock Convex envelopes of bivariate functions through the solution of {KKT}
  systems.
\newblock \emph{Journal of Global Optimization}, 72\penalty0 (2):\penalty0
  277--303, 2018.
\newblock \doi{10.1007/s10898-018-0626-1}.

\bibitem[Mahajan and Munson(2010)]{MahajanMunson2010}
A.~Mahajan and T.~Munson.
\newblock Exploiting second-order cone structure for global optimization.
\newblock Technical Report ANL/MCS-P1801-1010, Argonne National Laboratory,
  2010.

\bibitem[Maher et~al.(2017)Maher, Fischer, Gally, Gamrath, Gleixner, Gottwald,
  Hendel, Koch, L{\"u}bbecke, Miltenberger, M{\"u}ller, Pfetsch, Puchert,
  Rehfeldt, Schenker, Schwarz, Serrano, Shinano, Weninger, Witt, and
  Witzig]{SCIPoptsuite40}
S.~J. Maher, T.~Fischer, T.~Gally, G.~Gamrath, A.~Gleixner, R.~L. Gottwald,
  G.~Hendel, T.~Koch, M.~E. L{\"u}bbecke, M.~Miltenberger, B.~M{\"u}ller, M.~E.
  Pfetsch, C.~Puchert, D.~Rehfeldt, S.~Schenker, R.~Schwarz, F.~Serrano,
  Y.~Shinano, D.~Weninger, J.~T. Witt, and J.~Witzig.
\newblock {The SCIP Optimization Suite 4.0}.
\newblock ZIB Report 17-12, Zuse Institute Berlin, 2017.
\newblock \urnlink{nbn:de:0297-zib-62170}.

\bibitem[Margot(2010)]{Margot2010}
F.~Margot.
\newblock Symmetry in integer linear programming.
\newblock In M.~J{\"u}nger, T.~M. Liebling, D.~Naddef, G.~L. Nemhauser, W.~R.
  Pulleyblank, G.~Reinelt, G.~Rinaldi, and L.~A. Wolsey, editors, \emph{50
  Years of Integer Programming}, pages 647--686, Berlin, Heidelberg, 2010.
  Springer.
\newblock \doi{10.1007/978-3-540-68279-0\_17}.

\bibitem[McCormick(1976)]{McCormick1976}
G.~P. McCormick.
\newblock Computability of global solutions to factorable nonconvex programs:
  Part {I} -- convex underestimating problems.
\newblock \emph{Mathematical Programming}, 10\penalty0 (1):\penalty0 147--175,
  1976.
\newblock \doi{10.1007/bf01580665}.

\bibitem[MINLPLib()]{minlplib}
MINLPLib.
\newblock A library of mixed-integer and continuous nonlinear programming
  instances.
\newblock \url{https://www.minlplib.org}, 2022-10-14.

\bibitem[Misener and Floudas(2012)]{MisenerFloudas2012}
R.~Misener and C.~A. Floudas.
\newblock Global optimization of mixed-integer qua\-dra\-ti\-cally-constrained
  quadratic programs {(MIQCQP)} through piecewise-linear and edge-concave
  relaxations.
\newblock \emph{Mathematical Programming}, 136\penalty0 (1):\penalty0 155--182,
  2012.
\newblock \doi{10.1007/s10107-012-0555-6}.

\bibitem[Misener et~al.(2015)Misener, Smadbeck, and
  Floudas]{MisenerSmadbeckFloudas2015}
R.~Misener, J.~B. Smadbeck, and C.~A. Floudas.
\newblock Dynamically generated cutting planes for mixed-integer quadratically
  constrained quadratic programs and their incorporation into {GloMIQO} 2.
\newblock \emph{Optimization Methods and Software}, 30\penalty0 (1):\penalty0
  215--249, 2015.
\newblock \doi{10.1080/10556788.2014.916287}.

\bibitem[Moore(1966)]{Moore1966}
R.~E. Moore.
\newblock \emph{Interval Analysis}.
\newblock Englewood Cliffs, NJ: Prentice Hall, 1966.

\bibitem[M\"{u}ller et~al.(2020)M\"{u}ller, Serrano, and
  Gleixner]{MuellerSerranoGleixner2020}
B.~M\"{u}ller, F.~Serrano, and A.~Gleixner.
\newblock Using two-dimensional projections for stronger separation and
  propagation of bilinear terms.
\newblock \emph{{SIAM} Journal on Optimization}, 30\penalty0 (2):\penalty0
  1339--1365, 2020.
\newblock \doi{10.1137/19m1249825}.

\bibitem[Mu{\~{n}}oz and Serrano(2020)]{MunozSerrano2020}
G.~Mu{\~{n}}oz and F.~Serrano.
\newblock Maximal quadratic-free sets.
\newblock In D.~Bienstock and G.~Zambelli, editors, \emph{Integer Programming
  and Combinatorial Optimization}, pages 307--321, Cham, 2020. Springer.
\newblock \doi{10.1007/978-3-030-45771-6\_24}.

\bibitem[Müller et~al.(2018)Müller, Kuhlmann, and
  Vigerske]{MuellerKuhlmannVigerske2017}
B.~Müller, R.~Kuhlmann, and S.~Vigerske.
\newblock On the performance of {NLP} solvers within global {MINLP} solvers.
\newblock In \emph{Operations Research Proceedings 2017}, pages 633--639, Cham,
  2018. Springer International Publishing.
\newblock \doi{10.1007/978-3-319-89920-6\_84}.

\bibitem[Pfetsch(2022)]{marcscipversions}
M.~Pfetsch.
\newblock {SCIP: Past, Present, Future}.
\newblock Presentation at workshop \emph{Let's SCIP it!}, November 2022.
\newblock URL \url{https://scipopt.org/20years/slides/pfetsch.pdf}.

\bibitem[Pint\'er(2006)]{Pinter2006}
J.~D. Pint\'er, editor.
\newblock \emph{Global Optimization: Scientific and Engineering Case Studies},
  volume~85 of \emph{Nonconvex Optimization and Its Applications}.
\newblock Springer New York, NY, 2006.
\newblock \doi{10.1007/0-387-30927-6}.

\bibitem[Quesada and Grossmann(1993)]{QuGr93}
I.~Quesada and I.~E. Grossmann.
\newblock Global optimization algorithm for heat exchanger networks.
\newblock \emph{Industrial \& Engineering Chemistry Research}, 32\penalty0
  (3):\penalty0 487--499, 1993.
\newblock \doi{10.1021/ie00015a012}.

\bibitem[Sahinidis(1996)]{Sahinidis1996}
N.~Sahinidis.
\newblock {BARON}: {A} general purpose global optimization software package.
\newblock \emph{Journal of Global Optimization}, 8\penalty0 (2):\penalty0
  201--205, 1996.
\newblock \doi{10.1007/BF00138693}.

\bibitem[Schewe and Schmidt(2019)]{ScheweSchmidt2019}
L.~Schewe and M.~Schmidt.
\newblock Computing feasible points for binary {MINLPs} with {MPECs}.
\newblock \emph{Mathematical Programming Computation}, 11\penalty0
  (1):\penalty0 95--118, 2019.
\newblock \doi{10.1007/s12532-018-0141-x}.

\bibitem[Shinano(2021)]{ug}
Y.~Shinano.
\newblock {UG} -- {Ubiquity Generator framework} v1.0.0beta.
\newblock \url{https://ug.zib.de}, 2021.
\newblock \doi{10.12752/8521}.

\bibitem[Shinano et~al.(2016)Shinano, Achterberg, Berthold, Heinz, Koch, and
  Winkler]{ShinanoAchterbergBertholdHeinzKochWinkler2016}
Y.~Shinano, T.~Achterberg, T.~Berthold, S.~Heinz, T.~Koch, and M.~Winkler.
\newblock Solving open {MIP} instances with {ParaSCIP} on supercomputers using
  up to 80,000 cores.
\newblock In \emph{2016 IEEE International Parallel and Distributed Processing
  Symposium (IPDPS)}, pages 770--779, 2016.
\newblock \doi{10.1109/IPDPS.2016.56}.

\bibitem[Shinano et~al.(2018)Shinano, Heinz, Vigerske, and
  Winkler]{ShinanoHeinzVigerskeWinkler2018}
Y.~Shinano, S.~Heinz, S.~Vigerske, and M.~Winkler.
\newblock {FiberSCIP} -- a shared memory parallelization of {SCIP}.
\newblock \emph{INFORMS Journal on Computing}, 30\penalty0 (1):\penalty0
  11--30, 2018.
\newblock \doi{10.1287/ijoc.2017.0762}.

\bibitem[Smith and Pantelides(1999)]{SmithPantelides1999}
E.~Smith and C.~Pantelides.
\newblock A symbolic reformulation/spatial branch-and-bound algorithm for the
  global optimisation of nonconvex {MINLPs}.
\newblock \emph{Computers {\&} Chemical Engineering}, 23\penalty0
  (4-5):\penalty0 457--478, 1999.
\newblock \doi{10.1016/s0098-1354(98)00286-5}.

\bibitem[Smith et~al.(2013)Smith, Chinneck, and
  Aitken]{SmithChinneckAitken2013}
L.~Smith, J.~Chinneck, and V.~Aitken.
\newblock Improved constraint consensus methods for seeking feasibility in
  nonlinear programs.
\newblock \emph{Computational Optimization and Applications}, 54\penalty0
  (3):\penalty0 555--578, 2013.
\newblock \doi{10.1007/s10589-012-9473-z}.

\bibitem[Tardella(2004)]{Tardella2004}
F.~Tardella.
\newblock On the existence of polyhedral convex envelopes.
\newblock In C.~A. Floudas and P.~Pardalos, editors, \emph{Frontiers in Global
  Optimization}, pages 563--573, Boston, MA, 2004. Springer US.
\newblock \doi{10.1007/978-1-4613-0251-3\_30}.

\bibitem[Tawarmalani and Sahinidis(2005)]{TawarmalaniSahinidis2005}
M.~Tawarmalani and N.~V. Sahinidis.
\newblock A polyhedral branch-and-cut approach to global optimization.
\newblock \emph{Mathematical Programming}, 103\penalty0 (2):\penalty0 225--249,
  2005.
\newblock \doi{10.1007/s10107-005-0581-8}.

\bibitem[Trespalacios and Grossmann(2014)]{TrespalaciosGrossmann2014}
F.~Trespalacios and I.~Grossmann.
\newblock Review of mixed-integer nonlinear and generalized disjunctive
  programming methods.
\newblock \emph{Chemie Ingenieur Technik}, 86\penalty0 (7):\penalty0 991--1012,
  2014.
\newblock \doi{10.1002/cite.201400037}.

\bibitem[Tuy(1964)]{Tuy1964}
H.~Tuy.
\newblock Concave programming with linear constraints.
\newblock \emph{Doklady Akademii Nauk}, 159\penalty0 (1):\penalty0 32--35,
  1964.

\bibitem[Veinott(1967)]{Veinott1967}
A.~F. Veinott.
\newblock The supporting hyperplane method for unimodal programming.
\newblock \emph{Operations Research}, 15\penalty0 (1):\penalty0 147--152, 1967.
\newblock \doi{10.1287/opre.15.1.147}.

\bibitem[Vielma et~al.(2016)Vielma, Dunning, Huchette, and Lubin]{Vielma2016}
J.~P. Vielma, I.~Dunning, J.~Huchette, and M.~Lubin.
\newblock {Extended} formulations in mixed integer conic quadratic programming.
\newblock \emph{Mathematical Programming Computation}, 9\penalty0 (3):\penalty0
  369--418, 2016.
\newblock \doi{10.1007/s12532-016-0113-y}.

\bibitem[Vigerske(2013)]{Vigerske2013}
S.~Vigerske.
\newblock \emph{Decomposition of Multistage Stochastic Programs and a
  Constraint Integer Programming Approach to Mixed-Integer Nonlinear
  Programming}.
\newblock PhD thesis, Humboldt-Universit\"at zu Berlin, 2013.
\newblock \urnlink{nbn:de:kobv:11-100208240}.

\bibitem[Vigerske and Gleixner(2017)]{VigerskeGleixner2016}
S.~Vigerske and A.~Gleixner.
\newblock {SCIP}: global optimization of mixed-integer nonlinear programs in a
  branch-and-cut framework.
\newblock \emph{Optimization Methods and Software}, 33\penalty0 (3):\penalty0
  563--593, 2017.
\newblock \doi{10.1080/10556788.2017.1335312}.

\bibitem[W\"achter and Biegler(2006)]{WaechterBiegler2006}
A.~W\"achter and L.~Biegler.
\newblock On the implementation of a primal-dual interior point filter line
  search algorithm for large-scale nonlinear programming.
\newblock \emph{Mathematical Programming}, 106\penalty0 (1):\penalty0 25--57,
  2006.
\newblock \doi{10.1007/s10107-004-0559-y}.

\bibitem[Wegscheider(2019)]{Wegscheider2019}
F.~Wegscheider.
\newblock Exploiting symmetry in mixed-integer nonlinear programming.
\newblock Master's thesis, Zuse Institute Berlin, 2019.
\newblock \urnlink{nbn:de:0297-zib-77055}.

\bibitem[Xu et~al.(2022)Xu, D'Ambrosio, Liberti, and
  Vanier]{XuDAmbrosioLibertiVanier2022}
L.~Xu, C.~D'Ambrosio, L.~Liberti, and S.~H. Vanier.
\newblock On cutting planes for extended formulation of signomial programming.
\newblock Technical report, arXiv, 2022.
\newblock \doi{10.48550/ARXIV.2212.02857}.

\bibitem[Zamora and Grossmann(1998)]{ZamoraGrossmann1998}
J.~M. Zamora and I.~E. Grossmann.
\newblock Continuous global optimization of structured process systems models.
\newblock \emph{Computers and Chemical Engineering}, 22\penalty0 (12):\penalty0
  1749--1770, 1998.
\newblock \doi{10.1016/S0098-1354(98)00244-0}.

\end{thebibliography}

\appendix

\section{Test Set}
\label{sec:testsetdetail}

The following table provides details on the test set of 200 instances that was constructed by the selection process described in Section~\ref{sec:testset}.
For each instance, the number of variables ($n$), the number of discrete variables ($\vert\mathcal{I}\vert$), the number of constraints ($m+\tilde m$), the number of nonzeros in the Jacobian and objective function gradient (nz), and the number of nonzeros that correspond to nonlinear terms (nlnz) is given.

{
\scriptsize
\begin{longtable}{l|rrrrr}\toprule
instance & n & $\vert\mathcal{I}\vert$ & m+$\tilde m$ & nz & nlnz \\ \midrule
\endhead
alan & 8 & 4 & 7 & 23 & 3 \\
autocorr\_bern20-05 & 20 & 20 & 0 & 20 & 20 \\
autocorr\_bern35-04 & 35 & 35 & 0 & 35 & 35 \\
ball\_mk2\_10 & 10 & 10 & 1 & 20 & 10 \\
ball\_mk2\_30 & 30 & 30 & 1 & 60 & 30 \\
ball\_mk3\_10 & 10 & 10 & 1 & 20 & 10 \\
batch0812\_nc & 76 & 36 & 205 & 472 & 232 \\
batchs101006m & 278 & 129 & 1019 & 2865 & 49 \\
batchs121208m & 406 & 203 & 1511 & 4255 & 59 \\
bayes2\_20 & 86 & 0 & 77 & 615 & 440 \\
bayes2\_30 & 86 & 0 & 77 & 618 & 440 \\
blend029 & 102 & 36 & 213 & 542 & 64 \\
blend146 & 222 & 87 & 624 & 1721 & 256 \\
camshape100 & 199 & 0 & 200 & 696 & 299 \\
cardqp\_inlp & 50 & 50 & 1 & 100 & 50 \\
cardqp\_iqp & 50 & 50 & 1 & 100 & 50 \\
carton7 & 328 & 256 & 687 & 3979 & 678 \\
carton9 & 360 & 288 & 893 & 4917 & 758 \\
casctanks & 500 & 40 & 517 & 1605 & 514 \\
cecil\_13 & 840 & 180 & 898 & 2811 & 360 \\
celar6-sub0 & 640 & 640 & 16 & 1280 & 640 \\
chakra & 62 & 0 & 41 & 142 & 41 \\
chem & 11 & 0 & 4 & 36 & 11 \\
chenery & 43 & 0 & 38 & 132 & 56 \\
chimera\_k64maxcut-01 & 1101 & 1101 & 0 & 1101 & 1101 \\
chimera\_mis-01 & 2032 & 2032 & 0 & 2032 & 2032 \\
chp\_shorttermplan1a & 1008 & 144 & 2068 & 6118 & 576 \\
chp\_shorttermplan2a & 1584 & 240 & 3896 & 10160 & 1152 \\
chp\_shorttermplan2b & 1392 & 192 & 2552 & 7672 & 1440 \\
clay0204m & 52 & 32 & 90 & 284 & 64 \\
clay0205m & 80 & 50 & 135 & 430 & 80 \\
color\_lab3\_3x0 & 316 & 316 & 80 & 632 & 237 \\
crossdock\_15x7 & 210 & 210 & 44 & 630 & 210 \\
crossdock\_15x8 & 240 & 240 & 46 & 720 & 240 \\
crudeoil\_lee1\_07 & 749 & 56 & 1776 & 8124 & 896 \\
crudeoil\_pooling\_ct2 & 403 & 108 & 732 & 2523 & 140 \\
csched1 & 76 & 63 & 22 & 173 & 8 \\
csched1a & 28 & 15 & 22 & 77 & 7 \\
cvxnonsep\_psig20 & 20 & 10 & 0 & 20 & 20 \\
cvxnonsep\_psig30 & 30 & 15 & 0 & 30 & 30 \\
du-opt & 20 & 13 & 9 & 46 & 20 \\
du-opt5 & 20 & 13 & 9 & 46 & 20 \\
edgecross10-040 & 90 & 90 & 480 & 1530 & 90 \\
edgecross10-080 & 90 & 74 & 480 & 1528 & 88 \\
eg\_all\_s & 7 & 7 & 27 & 219 & 196 \\
eigena2 & 2500 & 0 & 1275 & 127500 & 127500 \\
elec50 & 150 & 0 & 50 & 300 & 300 \\
elf & 54 & 24 & 38 & 177 & 30 \\
eniplac & 141 & 24 & 189 & 510 & 48 \\
enpro56pb & 127 & 73 & 191 & 650 & 24 \\
ex1244 & 95 & 23 & 129 & 468 & 52 \\
ex1252a & 24 & 9 & 34 & 93 & 36 \\
faclay20h & 190 & 190 & 2280 & 7030 & 190 \\
faclay80 & 3160 & 3160 & 164320 & 496120 & 3160 \\
feedtray & 97 & 7 & 91 & 450 & 282 \\
fin2bb & 588 & 175 & 618 & 9413 & 42 \\
flay04m & 42 & 24 & 42 & 154 & 4 \\
flay05m & 62 & 40 & 65 & 242 & 5 \\
flay06m & 86 & 60 & 93 & 350 & 6 \\
fo7\_ar25\_1 & 112 & 42 & 269 & 1054 & 14 \\
fo7\_ar3\_1 & 112 & 42 & 269 & 1054 & 14 \\
forest & 236 & 73 & 309 & 1013 & 178 \\
gabriel01 & 215 & 72 & 467 & 1789 & 512 \\
gabriel02 & 261 & 71 & 597 & 2608 & 1024 \\
gasnet & 90 & 10 & 69 & 266 & 130 \\
gasprod\_sarawak16 & 1526 & 38 & 2252 & 6453 & 1088 \\
gastrans582\_cold13\_95 & 2186 & 250 & 3732 & 8538 & 2139 \\
gastrans582\_mild11 & 2186 & 250 & 3732 & 8538 & 2139 \\
gear & 4 & 4 & 0 & 4 & 4 \\
gear2 & 28 & 24 & 4 & 32 & 4 \\
gear4 & 6 & 4 & 1 & 8 & 4 \\
genpooling\_lee1 & 49 & 9 & 82 & 369 & 128 \\
genpooling\_lee2 & 53 & 9 & 92 & 453 & 192 \\
ghg\_1veh & 29 & 12 & 37 & 130 & 91 \\
gilbert & 1000 & 0 & 1 & 2000 & 2000 \\
graphpart\_2g-0066-0066 & 108 & 108 & 36 & 216 & 108 \\
graphpart\_clique-60 & 180 & 180 & 60 & 360 & 180 \\
gsg\_0001 & 77 & 0 & 111 & 368 & 44 \\
hadamard\_5 & 25 & 25 & 0 & 25 & 25 \\
heatexch\_spec1 & 56 & 12 & 64 & 224 & 32 \\
heatexch\_spec2 & 76 & 16 & 90 & 300 & 42 \\
hhfair & 29 & 0 & 25 & 80 & 21 \\
himmel16 & 18 & 0 & 21 & 96 & 84 \\
house & 8 & 0 & 8 & 25 & 9 \\
hs62 & 3 & 0 & 1 & 6 & 6 \\
hvb11 & 9817 & 9537 & 10251 & 36005 & 64 \\
hybriddynamic\_var & 81 & 10 & 100 & 286 & 61 \\
hybriddynamic\_varcc & 151 & 0 & 110 & 388 & 101 \\
hydroenergy1 & 288 & 96 & 428 & 1212 & 120 \\
ibs2 & 3010 & 1500 & 1821 & 13510 & 3000 \\
johnall & 194 & 190 & 192 & 957 & 573 \\
kall\_circles\_c6b & 17 & 0 & 53 & 148 & 86 \\
kall\_congruentcircles\_c72 & 17 & 0 & 59 & 160 & 86 \\
kissing2 & 772 & 0 & 10000 & 154400 & 154400 \\
kport20 & 101 & 40 & 27 & 189 & 116 \\
kriging\_peaks-red020 & 2 & 0 & 0 & 2 & 2 \\
kriging\_peaks-red100 & 2 & 0 & 0 & 2 & 2 \\
lop97icx & 986 & 899 & 87 & 1890 & 704 \\
mathopt5\_7 & 1 & 0 & 0 & 1 & 1 \\
mathopt5\_8 & 1 & 0 & 0 & 1 & 1 \\
maxcsp-geo50-20-d4-75-36 & 1000 & 1000 & 50 & 2000 & 1000 \\
meanvar-orl400\_05\_e\_7 & 2000 & 400 & 2003 & 7200 & 1600 \\
meanvar-orl400\_05\_e\_8 & 1600 & 400 & 1603 & 6400 & 800 \\
mhw4d & 5 & 0 & 3 & 13 & 10 \\
milinfract & 1000 & 500 & 501 & 502000 & 1000 \\
minlphi & 64 & 0 & 79 & 206 & 36 \\
multiplants\_mtg1a & 193 & 93 & 256 & 1972 & 95 \\
multiplants\_mtg2 & 229 & 112 & 306 & 2689 & 126 \\
nd\_netgen-3000-1-1-b-b-ns\_7 & 15000 & 3000 & 12155 & 48000 & 9000 \\
netmod\_kar1 & 456 & 136 & 666 & 1848 & 4 \\
netmod\_kar2 & 456 & 136 & 666 & 1848 & 4 \\
nous1 & 50 & 2 & 43 & 196 & 122 \\
nous2 & 50 & 2 & 43 & 196 & 122 \\
nvs02 & 8 & 5 & 3 & 19 & 16 \\
nvs06 & 2 & 2 & 0 & 2 & 2 \\
oil2 & 936 & 2 & 926 & 2214 & 440 \\
optmass & 30010 & 0 & 25005 & 80020 & 10006 \\
ortez & 87 & 18 & 74 & 268 & 54 \\
p\_ball\_10b\_5p\_3d\_m & 95 & 50 & 129 & 518 & 150 \\
p\_ball\_15b\_5p\_2d\_m & 105 & 75 & 139 & 523 & 150 \\
parabol5\_2\_3 & 40400 & 0 & 40200 & 240004 & 601 \\
parallel & 205 & 25 & 115 & 751 & 155 \\
pedigree\_ex485 & 485 & 426 & 296 & 1925 & 485 \\
pedigree\_ex485\_2 & 485 & 426 & 296 & 1925 & 485 \\
pointpack06 & 12 & 0 & 20 & 86 & 60 \\
pointpack08 & 16 & 0 & 35 & 155 & 112 \\
pooling\_epa1 & 214 & 30 & 340 & 1154 & 257 \\
pooling\_epa2 & 331 & 45 & 524 & 1913 & 554 \\
portfol\_buyin & 17 & 8 & 19 & 58 & 16 \\
portfol\_card & 17 & 8 & 20 & 66 & 16 \\
powerflow0014r & 118 & 0 & 197 & 652 & 461 \\
powerflow0057r & 440 & 0 & 725 & 2462 & 1795 \\
prob07 & 14 & 0 & 35 & 109 & 63 \\
process & 10 & 0 & 7 & 27 & 11 \\
procurement1mot & 784 & 60 & 749 & 2444 & 12 \\
procurement2mot & 796 & 60 & 761 & 2480 & 12 \\
product & 1553 & 107 & 1925 & 5555 & 264 \\
product2 & 2842 & 128 & 3125 & 8249 & 1056 \\
prolog & 20 & 0 & 22 & 128 & 14 \\
qp3 & 100 & 0 & 52 & 2747 & 100 \\
qspp\_0\_10\_0\_1\_10\_1 & 180 & 180 & 100 & 540 & 180 \\
qspp\_0\_11\_0\_1\_10\_1 & 220 & 220 & 121 & 660 & 220 \\
radar-2000-10-a-6\_lat\_7 & 10000 & 2000 & 8001 & 28000 & 6000 \\
radar-3000-10-a-8\_lat\_7 & 15000 & 3000 & 12001 & 42000 & 9000 \\
ravempb & 112 & 54 & 186 & 610 & 28 \\
risk2bpb & 463 & 14 & 580 & 2288 & 3 \\
routingdelay\_bigm & 1123 & 396 & 2977 & 7739 & 1827 \\
rsyn0815m & 205 & 79 & 347 & 909 & 11 \\
rsyn0815m03m & 705 & 282 & 1647 & 4120 & 33 \\
sfacloc2\_2\_95 & 186 & 39 & 239 & 595 & 76 \\
sfacloc2\_3\_90 & 291 & 75 & 496 & 1282 & 135 \\
sjup2 & 1696 & 8 & 17085 & 151716 & 88800 \\
slay06m & 102 & 60 & 135 & 462 & 12 \\
slay07m & 140 & 84 & 189 & 644 & 14 \\
smallinvDAXr1b010-011 & 30 & 30 & 3 & 120 & 30 \\
smallinvDAXr1b020-022 & 30 & 30 & 3 & 120 & 30 \\
sonet17v4 & 136 & 136 & 2057 & 6527 & 272 \\
sonet18v6 & 153 & 153 & 2466 & 7802 & 306 \\
sonetgr17 & 152 & 152 & 152 & 694 & 302 \\
spectra2 & 69 & 30 & 72 & 408 & 240 \\
sporttournament24 & 276 & 276 & 0 & 276 & 276 \\
sporttournament30 & 435 & 435 & 0 & 435 & 435 \\
sssd12-05persp & 95 & 75 & 52 & 305 & 45 \\
sssd18-06persp & 150 & 126 & 66 & 474 & 54 \\
st\_testgr1 & 10 & 10 & 5 & 51 & 10 \\
st\_testgr3 & 20 & 20 & 20 & 181 & 20 \\
steenbrf & 468 & 0 & 108 & 972 & 108 \\
stockcycle & 480 & 432 & 97 & 1008 & 48 \\
supplychainp1\_022020 & 2940 & 460 & 5300 & 15040 & 40 \\
supplychainp1\_030510 & 445 & 70 & 835 & 2330 & 15 \\
supplychainr1\_022020 & 1440 & 460 & 1840 & 7000 & 40 \\
supplychainr1\_030510 & 230 & 70 & 280 & 1005 & 15 \\
syn15m04m & 340 & 120 & 806 & 1986 & 44 \\
syn30m02m & 320 & 120 & 604 & 1502 & 40 \\
synheat & 56 & 12 & 64 & 224 & 28 \\
tanksize & 46 & 9 & 73 & 290 & 63 \\
telecomsp\_pacbell & 3570 & 3528 & 2940 & 121302 & 74088 \\
tln5 & 35 & 35 & 30 & 155 & 50 \\
tln7 & 63 & 63 & 42 & 287 & 98 \\
tls2 & 37 & 33 & 24 & 209 & 8 \\
tls4 & 105 & 89 & 64 & 613 & 32 \\
topopt-mbb\_60x40\_50 & 33600 & 2400 & 14363 & 259956 & 33600 \\
toroidal2g20\_5555 & 400 & 400 & 0 & 400 & 400 \\
toroidal3g7\_6666 & 343 & 343 & 0 & 343 & 343 \\
transswitch0009r & 69 & 9 & 103 & 346 & 255 \\
tricp & 169 & 0 & 190 & 1493 & 1140 \\
tspn08 & 44 & 28 & 18 & 136 & 60 \\
tspn15 & 135 & 105 & 34 & 502 & 165 \\
unitcommit1 & 960 & 720 & 5329 & 12404 & 240 \\
unitcommit2 & 960 & 720 & 5329 & 12404 & 480 \\
wager & 156 & 84 & 142 & 532 & 240 \\
waste & 2484 & 400 & 1991 & 9242 & 2736 \\
wastepaper3 & 52 & 27 & 30 & 177 & 108 \\
wastepaper4 & 76 & 44 & 38 & 274 & 176 \\
wastepaper6 & 136 & 90 & 54 & 528 & 360 \\
water4 & 195 & 126 & 137 & 756 & 46 \\
waternd1 & 74 & 20 & 83 & 301 & 114 \\
waterno2\_02 & 332 & 18 & 410 & 1088 & 202 \\
waterno2\_03 & 498 & 27 & 616 & 1635 & 303 \\
waterund01 & 40 & 0 & 38 & 152 & 78 \\
\bottomrule
\end{longtable}

}

\section{Detailed Computational Results}

The following tables show the outcome from running each solver on instances from the test set.
If an instance has been solved to optimality, the time spend is reported.
Note that due to differences in formulas for the relative gap in the various solvers, an instance may be accounted as solved even though the solver stopped at the time limit.
If a run has been flagged as failed, the reason for this decision is given: ``abort'' if the solver did not return with a result, ``nonopt'' if the reported upper or lower bound were not consistent with those given by MINLPLib, and ``infeas'' if the reported solution is not feasible with respect to the feasibility tolerance.
Otherwise, the relative gap at termination is reported, which is $\infty$ if no feasible solution or lower bound has been computed.
An exception here is BARON, where an instance is considered as solved if the solver only decided to not return a lower bound due to singularities in functions (see Section~\ref{sec:missingbounds}).
This is the case for instances \texttt{mhw4d} and \texttt{multiplants\_mtg2} and their permutations.
For each instance, a time or gap that is at most 10\% worse than the one from the best solver on this instance is printed in bold font.

\subsection{Serial Mode}
\label{sec:detailed_singlethread}

The following table shows the outcome from running each solver on the test set of 200 instances and their permutations in serial mode.

{
\scriptsize
\begin{longtable}{l@{}c|rrrr}
\toprule
instance & perm & BARON & Lindo API & Octeract & SCIP \\ \midrule
\endhead
alan & - & 0.33 & 0.89\% & \textbf{0.05} & 0.19\\
& 1 & 0.30 & 0.89\% & \textbf{0.05} & 0.14\\
& 2 & 0.19 & 0.89\% & 0.13 & \textbf{0.11}\\
& 3 & 0.19 & 0.89\% & \textbf{0.05} & 0.11\\
& 4 & 0.18 & 0.89\% & \textbf{0.05} & 0.13\\
autocorr\_bern20-05 & - & 7.29 & 100.43 & \textbf{4.36} & 16.34\\
& 1 & 5.18 & 90.50 & \textbf{3.73} & 13.26\\
& 2 & 6.19 & 86.68 & \textbf{4.58} & 23.15\\
\bottomrule
\end{longtable}

}

\subsection{Octeract Gap Limit}
\label{sec:detailed_octeractconvtol}

The following table shows the outcome from running Octeract in serial mode with both gap limits set to either $10^{-6}$ or $10^{-4}$ on the test set of 200 non-permuted instances.

{
\scriptsize
\begin{longtable}{l|rr}
\toprule
instance &  $10^{-6}$ & $10^{-4}$ \\ \midrule
\endhead
alan & \textbf{0.05} & \textbf{0.05}\\
autocorr\_bern20-05 & \textbf{4.36} & \textbf{4.39}\\
autocorr\_bern35-04 & \textbf{12.7\%} & \textbf{12.7\%}\\
ball\_mk2\_10 & \textbf{0.01} & \textbf{0.01}\\
ball\_mk2\_30 & \textbf{0.02} & 0.02\\
ball\_mk3\_10 & \textbf{0.04} & \textbf{0.04}\\
batch0812\_nc & \textbf{8.37} & \textbf{8.63}\\
batchs101006m & \textbf{4.84} & \textbf{4.92}\\
\bottomrule
\end{longtable}

}

\subsection{Parallel Mode}
\label{sec:detailed_multithread}

The following table shows the outcome from running BARON in serial and parallel mode.

{
\scriptsize
\begin{longtable}{l|rrrr}
\toprule
instance &  1 thread & 4 threads & 8 threads & 16 threads \\ \midrule
\endhead
alan & \textbf{0.33} & 0.38 & 0.37 & \textbf{0.36}\\
autocorr\_bern20-05 & 7.29 & 1.91 & \textbf{1.33} & \textbf{1.31}\\
autocorr\_bern35-04 & 29.91 & 8.44 & 5.71 & \textbf{4.66}\\
ball\_mk2\_10 & \textbf{0.07} & \textbf{0.07} & \textbf{0.07} & 0.20\\
ball\_mk2\_30 & 0.25 & \textbf{0.08} & \textbf{0.08} & \textbf{0.08}\\
ball\_mk3\_10 & \textbf{10.58} & \textbf{11.19} & \textbf{11.35} & 11.71\\
batch0812\_nc & 6.57 & \textbf{5.33} & 6.47 & \textbf{5.32}\\
batchs101006m & \textbf{6.91} & 13.82 & 9.41 & 11.02\\
\bottomrule
\end{longtable}

}

The following table shows the outcome from running Lindo API in serial and parallel mode.

{
\scriptsize
\begin{longtable}{l|rrrr}
\toprule
instance &  1 thread & 4 threads & 8 threads & 16 threads \\ \midrule
\endhead
alan & \textbf{0.89\%} & \textbf{0.89\%} & \textbf{0.89\%} & \textbf{0.89\%}\\
autocorr\_bern20-05 & 100.43 & 33.83 & \textbf{23.88} & \textbf{21.91}\\
autocorr\_bern35-04 & 13.5\% & 2655.83 & \textbf{1788.33} & 2461.81\\
ball\_mk2\_10 & 3.98 & 2.17 & \textbf{1.28} & 1.43\\
ball\_mk2\_30 & \textbf{100.0\%} & \textbf{100.0\%} & \textbf{100.0\%} & \textbf{100.0\%}\\
ball\_mk3\_10 & 3.58 & \textbf{1.79} & 2.05 & 2.17\\
batch0812\_nc & 25.46 & \textbf{18.73} & \textbf{17.51} & \textbf{17.68}\\
batchs101006m & \textbf{134.04} & abort & abort & abort\\
\bottomrule
\end{longtable}

}

The following table shows the outcome from running Octeract in serial and parallel mode.

{
\scriptsize
\begin{longtable}{l|rrrr}
\toprule
instance &  1 thread & 4 threads & 8 threads & 16 threads \\ \midrule
\endhead
alan & \textbf{0.05} & 0.06 & 0.08 & 0.09\\
autocorr\_bern20-05 & \textbf{4.36} & 4.91 & 4.88 & 7.71\\
autocorr\_bern35-04 & 12.7\% & 4064.07 & \textbf{1891.23} & 2241.67\\
ball\_mk2\_10 & \textbf{0.01} & 0.01 & 0.02 & 0.02\\
ball\_mk2\_30 & \textbf{0.02} & 0.02 & 0.02 & 0.03\\
ball\_mk3\_10 & \textbf{0.04} & \textbf{0.04} & 0.04 & 0.04\\
batch0812\_nc & 8.37 & 4.85 & \textbf{3.34} & 4.12\\
batchs101006m & \textbf{4.84} & 5.55 & 5.40 & 7.16\\
\bottomrule
\end{longtable}

}

The following table shows the outcome from running SCIP (1 thread) and FiberSCIP (4, 8, 16 threads).

{
\scriptsize
\begin{longtable}{l|rrrr}
\toprule
instance &  1 thread & 4 threads & 8 threads & 16 threads \\ \midrule
\endhead
alan & \textbf{0.19} & 1.05 & 1.06 & 1.13\\
autocorr\_bern20-05 & \textbf{16.34} & 33.08 & \textbf{17.10} & 18.19\\
autocorr\_bern35-04 & 107.46 & 158.09 & 101.10 & \textbf{64.17}\\
ball\_mk2\_10 & \textbf{0.05} & 1.05 & 1.06 & 1.12\\
ball\_mk2\_30 & \textbf{0.06} & 1.06 & 1.06 & 1.14\\
ball\_mk3\_10 & 0.00 & 0.00 & 0.00 & \textbf{0.00}\\
batch0812\_nc & 1.45 & \textbf{1.06} & \textbf{1.08} & \textbf{1.13}\\
batchs101006m & \textbf{7.63} & abort & abort & abort\\
\bottomrule
\end{longtable}

}

\end{document}